\documentclass[11pt]{article}
\usepackage{fullpage}
\usepackage{graphicx,amssymb,amsmath,theorem}
\usepackage{epsfig}
\usepackage{epsf}
\usepackage{times}
\usepackage{amsfonts}
\usepackage{latexsym}
\usepackage{bbm}     

\newtheorem{theorem}{Theorem}
\newtheorem{lemma}{Lemma}

\long\def\ignore#1{}
\def\cut{\ignore}
\newcommand{\later}[1]{{}}
\newcommand{\old}[1]{{}}

\newenvironment{proofof}[1]
      {\medskip\noindent{\bf Proof of #1:}\hspace{1mm}}
      {\hfill$\Box$\medskip}

\newcommand{\NN}{\mathbb{N}} 
\newcommand{\ZZ}{\mathbb{Z}} 
\newcommand{\RR}{\mathbb{R}} 

\def\reals{{\mathbb R}}

\def\A{{\cal A}}

\def\C{{\cal C}}
\def\L{{L}}

\def\T{{T}}

\def\bd{{\partial}}
\def\polylog{{\rm polylog}}
\def\eps{\varepsilon}
\def\etal{{\it et~al.}\,}
\def\uu{{\bf u}}
\def\vv{{\bf v}}
\def\xx{{\bf x}}
\def\kk{{\bf k}}

\def\marrow{{\marginpar[\hfill$\longrightarrow$]{$\longleftarrow$}}}

\newcommand\comm[1]{\typeout{Used \string\comm...}\marrow
\quad{\bf [[#1]]}\quad}

\title{Extremal problems on triangle areas in two and three dimensions}

\author{Adrian Dumitrescu\thanks{%
Department of Computer Science,
University of Wisconsin-Milwaukee, USA,
email: {\tt ad@cs.uwm.edu.}
Research partially supported by NSF CAREER grant CCF-0444188.}
\and
Micha Sharir\thanks{%
School of Computer Science, Tel Aviv University,
Tel Aviv, Israel, and Courant Institute, New York University, New
York, NY, USA,
email: {\tt michas@tau.ac.il.}
Research partially supported by NSF Grant CCR-05-14079,
by a grant from the U.S.-Israeli Binational Science Foundation,
by grant 155/05 from the Israel Science Fund, Israeli Academy of Sciences,
and by the Hermann Minkowski--MINERVA Center for Geometry at Tel Aviv
University.}
\and
Csaba D. T\'oth\thanks{%
Department of Mathematics,
Massachusetts Institute of Technology, Cambridge, MA, USA,
email: {\tt toth@math.mit.edu.}}}

\begin{document}

\maketitle
\thispagestyle{empty}

\begin{abstract}

The study of extremal problems on triangle areas was initiated in
a series of papers by Erd\H{o}s and Purdy in the early 1970s.
In this paper we present new results on such problems, concerning the
number of triangles of the same area that are spanned by finite point
sets in the plane and in 3-space, and the number of distinct areas 
determined by the triangles. 

In the plane, our main result is an $O(n^{44/19}) =O(n^{2.3158})$
upper bound on the number of unit-area triangles spanned by $n$
points, which is the first breakthrough improving the classical bound of
$O(n^{7/3})$ from 1992. We also make progress in a number of important
special cases: We show that (i) For points in convex position, there
exist $n$-element point sets that span $\Omega(n\log n)$ triangles of
unit area. (ii) The number of triangles of minimum (nonzero) area
determined by $n$ points is at most $\frac{2}{3}(n^2-n)$; there exist
$n$-element point sets (for arbitrarily large $n$)
that span $(6/\pi^2-o(1))n^2$ minimum-area
triangles.  (iii) The number of acute triangles of minimum area
determined by $n$ points is $O(n)$; this is asymptotically tight. (iv)
For $n$ points in convex position, the number of triangles of minimum
area is $O(n)$; this is asymptotically tight. (v) If no three points
are allowed to be collinear, there are $n$-element point sets that
span $\Omega(n\log n)$ minimum-area triangles
(in contrast to (ii), where collinearities are allowed and
a quadratic lower bound holds).

In 3-space we prove an $O(n^{17/7}\beta(n))= O(n^{2.4286})$ upper
bound on the number of unit-area triangles spanned by $n$ points,
where $\beta(n)$ is an extremely slowly growing function related to the
inverse Ackermann function. The best previous bound, $O(n^{8/3})$, is
an old result of Erd\H{o}s and Purdy from 1971. We further show,
for point sets in 3-space:
(i) The number of minimum nonzero area triangles is at most $n^2+O(n)$,
and this is worst-case optimal, up to a constant factor.
(ii) There are $n$-element point sets that span
$\Omega(n^{4/3})$ triangles of maximum area, all incident to a common
point. In any $n$-element point set, the maximum number of
maximum-area triangles incident to a common point is
$O(n^{4/3+\eps})$, for any $\eps>0$.
(iii) Every set of $n$ points, not all on a line, determines at least
$\Omega(n^{2/3}/\beta(n))$ triangles of distinct areas, which share a
common side.
\end{abstract}

\newpage
\setcounter{page}{1}

\section{Introduction}
\label{sec:intro}

Given $n$ points in the plane, consider the following equivalence
relation defined on the set of (nondegenerate) triangles  spanned by
the points: two triangles are {\em equivalent} if they have the
same area.
Extremal problems typically ask for the maximum cardinality of an
equivalence class, and for the minimum number of distinct equivalence
classes, in a variety of cases. A classical example is when we call
two segments spanned by the given points equivalent if they
have the same length. Bounding the maximum size of an equivalence class
is the famous {\em repeated distances} problem
\cite{bmp-05,e-46,sst-84,s-97}, and bounding the minimum number of distinct
classes is the equally famous
{\em distinct distances} problem \cite{bmp-05,e-46,kt-04,st-01,s-97,t-03}.
In this paper, we make progress on several old extremal problems on
triangle areas in two and in three dimensions. We also study some
new and interesting variants never considered before. Our proof
techniques draw from a broad range of combinatorial tools such as
the Szemer\'edi-Trotter theorem on point-line incidences~\cite{st-83},
the Crossing Lemma~\cite{acns-82,l-84}, incidences between curves and
points and tangencies between curves and lines, extremal graph
theory~\cite{kst-54}, quasi-planar graphs~\cite{aapps-97},
Minkowski-type constructions, repeated distances on the
sphere~\cite{pa-95}, the partition technique of Clarkson~\etal~\cite{ceg-90},
various charging schemes, etc.

In 1967, A.~Oppenheim (see \cite{ep-95}) asked the following question:
Given $n$ points in the plane and $A>0$, how many triangles spanned
by the points can have area $A$? By applying an affine transformation,
one may assume $A=1$ and count the triangles of {\em unit} area.
Erd\H{o}s and Purdy~\cite{ep-71} showed that a $\sqrt{\log n}\times
(n/\sqrt{\log n})$ section of the integer lattice determines
$\Omega(n^2 \log\log{n})$ triangles of the same area. They also showed
that the maximum number of such triangles is at most $O(n^{5/2})$. In
1992, Pach and Sharir~\cite{ps-92} improved the exponent and obtained
an $O(n^{7/3})$ upper bound using the Szemer\'edi-Trotter
theorem~\cite{st-83} on the number of point-line incidences. We
further improve the upper bound by estimating the number of incidences
between the points and a 4-parameter family of quadratic curves. We
show that $n$ points in the plane determine at most
$O(n^{44/19})=O(n^{2.3158})$ unit-area triangles. We also consider the
case of points in convex position, for which we construct $n$-element
point sets that span $\Omega(n\log n)$ triangles of unit area.

Bra\ss, Rote, and Swanepoel~\cite{brs-01} showed that $n$ points in
the plane determine at most $O(n^2)$ minimum-area triangles, and they
pointed out that this bound is asymptotically tight. We introduce a
simple charging scheme to first bring the upper bound down to
$n^2-n$ and then further to $\frac{2}{3}(n^2-n)$. Our charging scheme
is also instrumental in showing that a $\sqrt{n} \times \sqrt{n}$
section of the integer lattice spans $(6/\pi^2-o(1))n^2$ triangles of
minimum area.  In the lower bound constructions, there are many collinear
triples and most of the minimum-area triangles are obtuse. We show
that there are at most $O(n)$ {\em acute} triangles of minimum
(nonzero) area, for any $n$-element point set.
Also, we show that $n$ points in (strictly) 
convex position determine at most
$O(n)$ minimum-area triangles---these bounds are best possible apart from
the constant factors. If no three points are allowed to be collinear,
we construct $n$-element point sets that span $\Omega(n\log{n})$
triangles of minimum area.

Next we address analogous questions for triangles in 3-space.
The number of triangles with some extremal property
might go up (significantly) when one moves up one dimension.
For instance, Bra\ss, Rote, and Swanepoel~\cite{brs-01} have shown that
the number of maximum area triangles in the plane is at most $n$
(which is tight). In 3-space we show that this number is at
least $\Omega(n^{4/3})$ in the worst case.
In contrast, for minimum-area triangles, we
prove that the quadratic upper bound from the planar case remains in
effect for 3-space, with a different constant of proportionality.

As mentioned earlier, Erd\H{o}s and Purdy~\cite{ep-71} showed that a
suitable $n$-element section of the integer lattice determines
$\Omega(n^2 \log\log{n})$ triangles of the same area. Clearly, this
bound is also valid in 3-space. In the same paper, via a forbidden
graph argument applied to the incidence graph between points and
cylinders whose axes pass through the origin, Erd\H{o}s and Purdy
deduced an $O(n^{5/3})$ upper bound on the number of unit-area triangles
incident to a common point, and thereby an $O(n^{8/3})$ upper bound on
the number of unit-area triangles determined by $n$ points in 3-space.
Here, applying a careful (and somewhat involved) analysis of the
structure of point-cylinder incidences in $\RR^3$, we
prove a new upper bound of $O(n^{17/7}\beta(n))=O(n^{2.4286})$,
for $\beta(n) = \exp(\alpha(n)^{O(1)})$,
where $\alpha(n)$ is the
extremely slowly growing inverse Ackermann function.

It is conjectured~\cite{bmp-05,bp-79,eps-82} that $n$ points in
$\RR^3$, not all on a line, determine at least $\lfloor(n-1)/2\rfloor$
distinct triangle areas. This bound has recently been established in
the plane \cite{p-07}, but the question is still wide open in $\RR^3$. 
It is attained by $n$ equally spaced points distributed evenly on two
parallel lines (which is in fact a planar construction).
We obtain a first result on this question and show that
$n$ points in $\RR^3$, not all on a line, determine at least
$n^{2/3}\exp(-\alpha(n)^{O(1)})=\Omega(n^{.666})$ triangles of
distinct areas.  Moreover, all these triangles share a common side.

\section{Unit-area triangles in the plane}
\label{sec:unit2}

\noindent {\bf The general case.} We establish a new upper bound on the
maximum number of unit-area triangles determined by $n$ points the
plane. 

\begin{theorem}\label{thm:unit2}
The number of unit-area triangles spanned by $n$ points in the plane
is $O(n^{2+6/19})= O(n^{2.3158})$.
\end{theorem}
\begin{proof}
Let $S$ be a set of $n$ points in the plane. Consider a triangle
$\Delta{abc}$ spanned by $S$. We call the three lines containing the
three sides of $\Delta{abc}$, {\em base lines} of $\Delta$, and the
three lines parallel to the base lines and incident to the third
vertex, {\em top lines} of $\Delta$.

\cut{For any integer $k\in \NN$, let $L_k$ denote the set of $k$-rich
  lines (that is, the lines containing at least $k$ points of $S$). By
  the Szemer\'edi-Trotter theorem~\cite{st-83}, we have
  $|L_k|=O(n^2/k^3+n/k)$. In particular, $|L_k|=O(n^2/k^3)$ for $k\leq
  \sqrt{n}$ and $|L_k| = O(n/k)$ for $k\geq \sqrt{n}$.}
For a parameter $k$, $1\leq k\leq \sqrt{n}$, to be optimized
later, we partition the set of unit-area triangles as follows.

\smallskip
\noindent $\bullet$ $U_1$ denotes the set of unit-area triangles where
one of the top lines is incident to fewer than $k$ points of $S$.

\smallskip
\noindent $\bullet$ $U_2$ denotes the set of unit-area triangles where
all three top lines are {\em $k$-rich} (i.e., each contains at
least $k$\\
\indent points of $S$).

\smallskip
\noindent We derive different upper bounds for each of these
types of unit-area triangles.

\paragraph{Bound for  $|U_1|$.}
For any two distinct points, $a,b\in \RR^2$, let $\ell_{ab}$ denote
the line through $a$ and $b$. The points $c$ for which the triangle
$\Delta{abc}$ has unit area lie on two lines
$\ell^-_{ab},\ell^+_{ab}$ parallel to $\ell_{ab}$ at distances
$2/|ab|$ on either side of $\ell_{ab}$. The ${n\choose 2}$ segments
determined by $S$ generate at most $2{n\choose 2}$ such lines
(counted with multiplicity).
If $\Delta abc \in U_1$ and its top line incident to the fewest
points of $S$ is $\ell_{ab}'\in \{\ell^-_{ab},\ell^+_{ab}\}$, then
$\ell_{ab}'$ is incident to at most $k$ points, so the segment $ab$ is
the base of at most $k$ triangles $\Delta{abc}\in U_1$ (with
$c\in\ell_{ab}'$).  This gives the upper bound
$$|U_1|\leq  2{n\choose 2} \cdot k = O(n^2 k).$$

\paragraph{Bound for $|U_2|$.}
Let $L$ be the set of $k$-rich lines, and let $m=|L|$. By the
Szemer\'edi-Trotter theorem~\cite{st-83}, we have $m=O(n^2/k^3)$ for
any $k\leq \sqrt{n}$. Furthermore, the cardinality of the set
$I(S,L)$ of point-line incidences between $S$ and $L$ is
$|I(S,L)|=O(n^2/k^2)$.

\old{
If $\Delta{abc}\in U_2$, then its top lines
$\ell_{ab}',\ell_{bc}',\ell_{ac}'\in L_k$ form a triangle of area 4,
in which the points $a,b,c$ are the midpoints of the three sides. We
give an upper bound on the number of 4-tuples
$(\ell_1,\ell_2,\ell_3,a)\in L_k\times S$ where the lines
$\ell_1,\ell_2,\ell_3$ form a triangle of area 4 and $a$ is the
midpoint of a side. Since every triangle in $U_2$ corresponds to
three 4-tuples of this form, we obtain an upper bound on $|U_2$.
} 

For any pair of nonparallel lines $\ell_1,\ell_2\in L$, let
$\gamma(\ell_1,\ell_2)$ denote the locus of points $p\in \RR^2$,
$p\not\in\ell_1\cup \ell_2$, such that the parallelogram that has a
vertex at $p$ and two sides along $\ell_1$ and $\ell_2$,
respectively, has area 2. The set $\gamma(\ell_1,\ell_2)$ consists
of two hyperbolas with $\ell_1$ and $\ell_2$ as asymptotes.
See Figure~\ref{fig:hyper}. For instance, if
$\ell_1: y=0$ and $\ell_2: y= ax$, then $\gamma(\ell_1,\ell_2) =
\{(x,y)\in \RR^2:xy=y^2/a+2\}\cup \{(x,y)\in \RR^2:xy=y^2/a-2\}$.
Any two nonparallel lines uniquely determine two such hyperbolas.
Let $\Gamma$ denote the set of these hyperbolas. Note that
$|\Gamma|=O(m^2)$. The family of such hyperbolas for all pairs of
nonparallel lines form a 4-parameter family of quadratic curves
(where the parameters are the coefficients of the defining lines).

For any triangle $\Delta{abc}\in U_2$, any pair of its top lines,
say, $\ell_{ab}'$ and $\ell_{ac}'$, determine a hyperbola passing
through $a$, which is incident to the third top line $\ell_{bc}'$;
furthermore $\ell_{bc}'$ is tangent\footnote{%
  For a quick proof, let $\uu$ (resp., $\vv$) be a unit vector
  along $\ell'_{ac}$ (resp., $\ell'_{ab}$). The point $a$ can be
  parametrized as $\xx= t\uu + \frac{\kappa}{t}\vv$, where
  $\kappa=2/\sin\theta$, and $\theta$ is the angle between
  $\ell'_{ac}$ and $\ell'_{ab}$. Hence the tangent to the
  hyperbola at $a$ is given by
  $\dot{\xx} = \uu-\frac{\kappa}{t^2}\vv \,\|\, t\uu-\frac{\kappa}{t}\vv =
  \vec{cb}$.}
to the hyperbola at $a$. See Figure~\ref{fig:hyper}. Any hyperbola in
this 4-parameter family is uniquely determined by two incident
points and the two respective tangent lines at those points.

\begin{figure}[htbp]
  \begin{center}
\epsfig{file=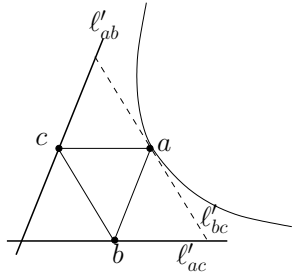,width=3.9cm,clip=}
 \caption{\small One of the hyperbolas defined by the triangle
 $\Delta abc$.
\label{fig:hyper}}
\end{center}
 \end{figure}

\old{If three lines $\ell_1,\ell_2,\ell_3\in L_k$ form a triangle
$\Delta$ of area 4, then $\ell_3$ is tangent to
$\gamma(\ell_1,\ell_3)$, and the intersection point
$\gamma(\ell_1,\ell_2)\cap \ell_3$ is the midpoint of the side of
$\Delta$ along $\ell_3$. Any hyperbola in the 4-parameter family
$\Gamma$ is uniquely determined by two incident points and two
tangent lines at those points.
} 

We define a topological graph $G$ as follows. For each point
$p\in S$, which is incident to $d_p$ lines of $L$, we create $2d_p$
vertices in $G$, as follows (refer to Figure~\ref{fig:graph}).
Draw a circle
$C_\eps(p)$ centered at $p$ with a sufficiently small radius
$\eps>0$, and place a vertex at every intersection point of the
circle $C_\eps(p)$ with the $d_p$ lines incident to $p$. The number
of vertices is $v_G=2|I(S,L)|=O(n^2/k^2)$. Next, we define the edges
of $G$. For each connected branch $\gamma$ of every hyperbola in
$\Gamma$, consider the set $S(\gamma)$ of points $p\in S$ that are
(i) incident to $\gamma$ and (ii) some line of $L$ is tangent to
$\gamma$ at $p$. For any two consecutive points $p,q\in S(\gamma)$,
draw an edge along $\gamma$ between the two vertices of $G$ that
(i) correspond to the incidences $(p,\ell_p)$ and $(q,\ell_q)$, where
$\ell_p$ and $\ell_q$ are the tangents of $\gamma$ at $p$ and $q$,
respectively, and (ii) are closest to each other along $\gamma$.
Specifically, the edge follows $\gamma$ between the
circles $C_{2\eps}(p)$ and $C_{2\eps}(q)$ and follows straight line
segments in the interiors of those circles.
Choose $\eps>0$ sufficiently small so that the circles $C_{2\eps}(p)$
have disjoint interiors and
the portions of the hyperbolas in the interiors of the circles
$C_{2\eps}(p)$, for every $p\in S$, meet at $p$ only. This
guarantees that the edges of $G$ cross only at intersection points
of the hyperbolas. The graph $G$ is {\em simple} because two points and
two tangent lines uniquely determine a hyperbola in $\Gamma$. The
number of edges is at least $3|U_2|-2m^2$, since every triangle in
$U_2$ corresponds to three point-hyperbola incidences in $I(S,\Gamma)$
(satisfying the additional condition of tangency with the respective
top lines);
and along each of the $2m^2$ hyperbola branches, each of its
incidences with the points of $S$ (of the special kind under
consideration), except for one, contributes one edge to $G$.
\begin{figure}[htbp]
  \begin{center}
\epsfig{file=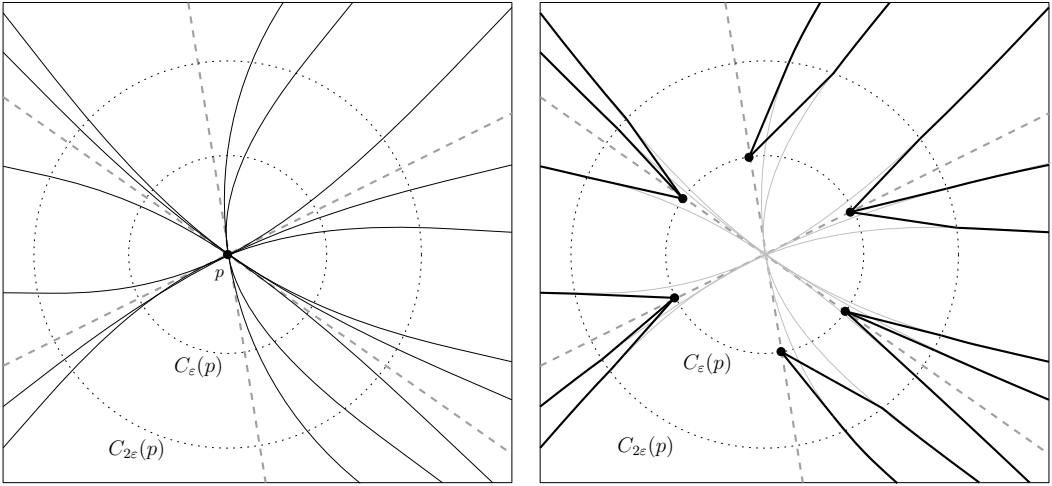,width=14cm,clip=}
 \caption{\small On the left: a point $p\in S$ incident to three lines
of $L$ (dashed) and 8 hyperbolas, each tangent to one of those
lines. On the right: the 6 vertices of $G$ corresponding to the 3
point-line incidences at $p$, and the drawings of the edges along the
hyperbolas.\label{fig:graph}}
\end{center}
 \end{figure}
Thus $G$ is a simple topological graph with $v_G=2I(S,L)=O(n^2/k^2)$
vertices and $e_G \geq 3|U_2|-2m^2$ edges. Since in this drawing of
$G$, every crossing is an intersection of two hyperbolas, the
crossing number of $G$ is upper bounded by ${\rm cr}(G)= O(|\Gamma|^2) =
O(m^4)$. We can also bound the crossing number of $G$ from below via
the Crossing Lemma of Ajtai~\etal~\cite{acns-82} and
Leighton~\cite{l-84}. It follows that
$$\Omega\left( \frac{e_G^3}{v_G^2}\right) -4v_G \leq {\rm cr}(G) \leq
O(m^4).$$
Rearranging this chain of inequalities, we obtain $e_G^3 = O(m^4v_G^2
+ v_G^3)$, or $e_G= O(m^{4/3}v_G^{2/3} + v_G)$. Comparing this bound
with our lower bound $e_G\ge 3|U_2|-2m^2$, we have $|U_2|=
O(m^{4/3}v_G^{2/3} + v_G+ m^2)$. Hence, for $k\leq \sqrt{n}$, we
have
$$|U_2|= O\left( \left(\frac{n^2}{k^3}\right)^{4/3}
  \left(\frac{n^2}{k^2}\right)^{2/3} +\frac{n^2}{k^2} +
  \left(\frac{n^2}{k^3} \right)^2\right) =
O\left(\frac{n^4}{k^{16/3}}+\frac{n^2}{k^2}\right)=
O\left(\frac{n^4}{k^{16/3}}\right).$$
The total number of unit-area triangles is $|U_1|+|U_2| = O(n^2k+
n^4/k^{16/3})$. This expression is minimized for
$k=n^{6/19}$, and we get $|U_1|+|U_2| = O(n^{44/19})$.
\end{proof}

\subsection{Convex position}

The construction of Erd\H{o}s and Purdy~\cite{ep-71} with many triangles
of the same area, the $\sqrt{\log n}\times (n/\sqrt{\log n})$ section of
the integer lattice, also contains many collinear triples. Here we consider
the unit-area triangle problem in the special case of point sets in
strictly convex position, so no three points are collinear.
We show that $n$ points in convex position in the plane can determine a
superlinear number of unit-area triangles. On the other hand, we do not know
of any subquadratic upper bound.

\begin{theorem} \label{thm:convex-unit}
For all $n \geq 3$, there exist $n$-element point sets in convex
position in the plane that span $\Omega(n\log n)$ unit-area
triangles.
\end{theorem}
\begin{proof}
We recursively construct a set $S_i$ of $n_i=3^i$ points on the
unit circle that determine $t_i=i3^{i-1}$ unit-area triangles, for
$i=1,2,\ldots$.  Take a circle $C$ of unit radius centered
at the origin $o$. We start with a set $S_1$ of 3 points along the
circle forming a unit-area triangle, so we have $n_1=3$ points and
$t_1=1$ unit-area triangles. In each step, we triple the number of
points, i.e., $n_{i+1}=3n_i$, and create new unit-area triangles, so that
$t_{i+1}=3t_i+n_i$. This implies $n_i=3^i$, and $t_i= i 3^{i-1}$,
yielding the desired lower bound.
The $i$-th step, $i \geq 2$, goes as follows. Choose a
generic angle value $\alpha_i$, close to $\pi/2$, say,
and let $\beta_i$ be the angle such that the
three unit vectors at direction $0$, $\alpha_i$, and $\beta_i$ from
the origin determine a unit-area triangle, which we denote by $D_i$
(note that $\beta_i$ lies in the third quadrant).
Rotate $D_i$ around the origin to each position where its $0$ vertex
coincides with one of the $n_i$ points of $S_i$, and add the other
two vertices of $D_i$ in these positions
to the point set. (With appropriate choices of $S_1$ and the angles
$\alpha_i$, $\beta_i$, one can guarantee that no two points of any
$S_i$ coincide.) For each point of $S_i$, we
added two new points, so $n_{i+1}=3n_i$. Also, we have $n_i$ new unit-area
triangles from rotated copies of $D_i$; and each of the $t_i$
previous triangles have now two new copies rotated by $\alpha_i$ and
$\beta_i$. This gives $t_{i+1}=3t_i+n_i$.
\end{proof}

\section{Minimum-area triangles in the plane}
\label{sec:minimum2}

\paragraph{The general case.} We first present a simple but
effective charging scheme that gives an upper bound of $n^2-n$ on the
number of minimum (nonzero) area triangles spanned by $n$ points in
the plane (Lemma \ref{lem:min-area}).
This technique yields a very short proof of the minimum area result
from \cite{brs-01}, with a much better constant of proportionality.
Moreover, its higher-dimensional variants lead to asymptotically tight
bounds on the maximum number of minimum-volume $k$-dimensional
simplices in $\RR^d$, for any $1 \leq k \leq d$
(see Section \ref{sec:minimum3} for the case $k=2,d=3$, and
\cite{dt-07b} for the case $k=3,d=3$; the generalization to arbitrary
$1\leq k \leq d$ will be presented in the journal version of \cite{dt-07b}).

\begin{lemma} \label{lem:min-area}
 The number of triangles of minimum (nonzero)
 area spanned by $n$ points in the plane is at most $n^2-n$.
\end{lemma}
\begin{proof}
Consider a set $S$ of $n$ points in the plane. Assign every triangle
of minimum area to one of its longest sides. For a segment $ab$,
with $a,b\in S$, let $R_{ab}^+$ and $R_{ab}^-$ denote the two
rectangles of extents $|ab|$ and $2/|ab|$ with $ab$ as a common side.
If a minimum-area triangle $\Delta{abc}$ is assigned to $ab$, then
$c$ must lie in the relative interior of the side parallel to $ab$ in either
$R_{ab}^+$ or $R_{ab}^-$. If there were two points, $c_1$ and $c_2$,
on one of these sides, then the area of $\Delta{ac_1c_2}$ would be smaller
than that of $\Delta{abc}$, a contradiction. Therefore, at most two
triangles are assigned to each of the ${n \choose 2}$ segments
(at most one on each side of the segments), and so there are at
most $n^2-n$ minimum-area triangles.
\end{proof}

\smallskip
We now refine our analysis and establish a $\frac{2}{3}(n^2-n)$ upper bound, which
leaves only a small gap from our lower bound
$(\frac{6}{\pi^2}-o(1))n^2$; both bounds are presented in
Theorem~\ref{thm:min-area} below. Let us point out again
that here we allow collinear triples of points. The maximum number of
collinear triples is clearly ${n\choose 3}=\Theta(n^3)$.
The bounds below, however, consider only nondegenerate triangles of
{\em positive} areas.


\begin{theorem} \label{thm:min-area}
 The number of triangles of minimum (nonzero)
 area spanned by $n$ points in the plane is at most
 $\frac{2}{3}(n^2-n)$.
The points in the $\lfloor\sqrt{n}\rfloor\times
\lfloor\sqrt{n}\rfloor$ integer grid span
$(\frac{6}{\pi^2}-o(1))n^2 \gtrapprox .6079 n^2$ minimum-area
triangles.
\end{theorem}
\begin{proof}
We start with the upper bound.
Consider a set $S$ of $n$ points in the plane, and let $L$ be the
set of connecting lines determined by $S$. Assume, without loss of
generality, that none of the lines in $L$ is vertical.
Let $T$ be the set of minimum (nonzero) area triangles spanned by $S$,
and put $t=|T|$.
There are $3t$ pairs $(ab,c)$ where $\Delta{abc}\in T$, and we may
assume, without loss of generality, that for at least half of these
pairs (i.e., for at least $\frac{3}{2}t$ pairs) $\Delta{abc}$ lies above
the line spanned by $a$ and $b$.

For each line $\ell\in L$, let $\ell'$ denote the line parallel to
$\ell$, lying above $\ell$, passing through some point(s) of $S$, and
closest to $\ell$ among these lines. Clearly, if $c\in S$
generates with $a,b\in\ell$ a minimum-area triangle which lies above
$ab$ then (i) $a$ and $b$ are a closest pair among the pairs of points
in $\ell\cap S$, and (ii) $c\in\ell_{ab}'$
(the converse does not necessarily hold).

Now fix a line $\ell\in L$; set $k_1=|\ell\cap S|\ge 2$, and
$k_2=|\ell'\cap S|\ge 1$, where $\ell'$ is as defined above.
The number of minimum-area triangles determined by a pair of points in
$\ell$ and lying above $\ell$ is at most $(k_1-1)k_2$. We have
\begin{equation} \label{eq:k12}
{k_1\choose 2} + {k_2\choose 2} \ge (k_1-1)k_2 .
\end{equation}
Indeed, multiplying by $2$ and subtracting the right-hand side from
the left-hand side gives
$$
k_1^2 -k_1 + k_2^2 -k_2 - 2k_1k_2 + 2k_2 =
(k_1-k_2)^2 - (k_1-k_2) \ge 0,
$$
which holds for any $k_1,k_2 \in \ZZ$.

We now sum (\ref{eq:k12}) over all lines $\ell\in L$. The sum of the
terms ${k_1\choose 2}$ is ${n\choose 2}$, and the sum of the terms
${k_2\choose 2}$ is at most ${n\choose 2}$, because a line $\lambda\in L$
spanned by at least two points of $S$ can arise as the line $\ell'$ for at most
one line $\ell\in L$. Hence we obtain
$$
\frac{3}{2}t \le \sum_{\ell\in L} (k_1-1)k_2 \le 2{n\choose 2} =
n(n-1) ,
$$
thus $t\le \frac23(n^2-n)$, as asserted.

\smallskip
We now prove the lower bound. Consider the set $S$ of points in the
$\lfloor \sqrt{n}\rfloor \times \lfloor\sqrt{n}\rfloor$ section of the
integer lattice. Clearly $|S| \leq n$. The minimum nonzero area of
triangles in $S$ is $1/2$ (by Pick's theorem). Recall that the
charging scheme used in the proof of Lemma~\ref{lem:min-area}
assigns each triangle of minimum area to one of its
longest sides, which is necessarily a {\em visibility segment}
(a segment not containing any point of $S$ in its relative interior).
We show that every visibility segment $ab$ which is not
axis-parallel is assigned to exactly two triangles of minimum area.

Draw parallel lines to $ab$ through all points of the integer
lattice. Every line parallel to $ab$ and incident to a point of $S$
contains equally spaced points of the (infinite) integer lattice.
The distance between consecutive points along each line is exactly
$|ab|$.
This implies that each of the two lines parallel to $ab$ and closest
to it contains a lattice point on the side of the respective rectangle
$R_{ab}^-$ or $R_{ab}^+$,
opposite to $ab$, and this lattice point is in $S$.
Finally, observe that there are no empty acute triangles in the integer
lattice. It follows that our charging scheme uniquely assigns empty
triangles to visibility segments. An illustration is provided in
Figure~\ref{fig:grid}.

\begin{figure}[htbp]
  \begin{center}
\epsfig{file=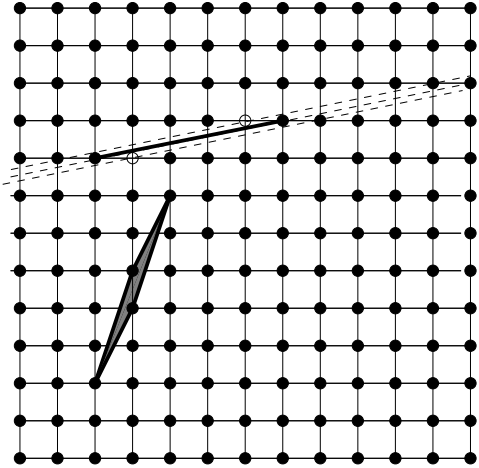,width=2.5in,clip=}
 \caption{\small In an integer lattice section, every visibility segment
 which is not axis-parallel is the longest side of two triangles of
 minimum area.}
\end{center}
\label{fig:grid}
\end{figure}

A non-axis-parallel segment $ab$ is a visibility segment if
and only if the coordinates of the vector $\overrightarrow{ab}$ are
relatively prime. It is well known that $6/\pi^2$ is the limit of the
probability that a pair of integers $(i,j)$ with $ 1 \leq i,j \leq m$
are relatively prime, as $m$ tends to infinity \cite{v-54}.
Hence, a fraction of about
$6/\pi^2$ of the ${|S| \choose 2} \leq {n \choose 2}$ segments
spanned by $S$ are visibility segments which are not axis-parallel.
Each of these $(\frac{6}{\pi^2}-o(1)){n\choose 2}$ segments corresponds
to two unique triangles of minimum area, so $S$ determines at least
$(\frac{6}{\pi^2}-o(1))n^2$ minimum-area triangles.
\end{proof}

\subsection{Special cases}

In this subsection we consider some new variants of the minimum-area
triangle problem for the two special cases (i) where no three points are
collinear, and (ii) where the points are in convex position.
We also show that the maximum number of {\em acute} triangles
of minimum area, for any point set, is only linear.

\paragraph{Acute triangles.}
We have seen that $n$ points in an integer grid may span
$\Omega(n^2)$ triangles of minimum area. However, in that
construction, all these triangles are obtuse (or right-angled).
Here we prove that for any $n$-element point set in the plane, the
number of {\em acute} triangles of minimum area is only linear.
    This bound is attained in the following simple example.
    Take two groups of about $n/2$ equally spaced points on two parallel
    lines: the first group consist of the points
    $(i,0)$, for $i=0,\ldots,\lceil n/2 \rceil -1$, and the second
    group of the
    points $(i+1/2,\sqrt{3}/2)$, for $i=0,\ldots,\lfloor n/2 \rfloor -1$.
    This point set determines $n-2$ acute triangles of minimum area.

\begin{figure}[htbp]
\centerline{\epsfxsize=3in \epsffile{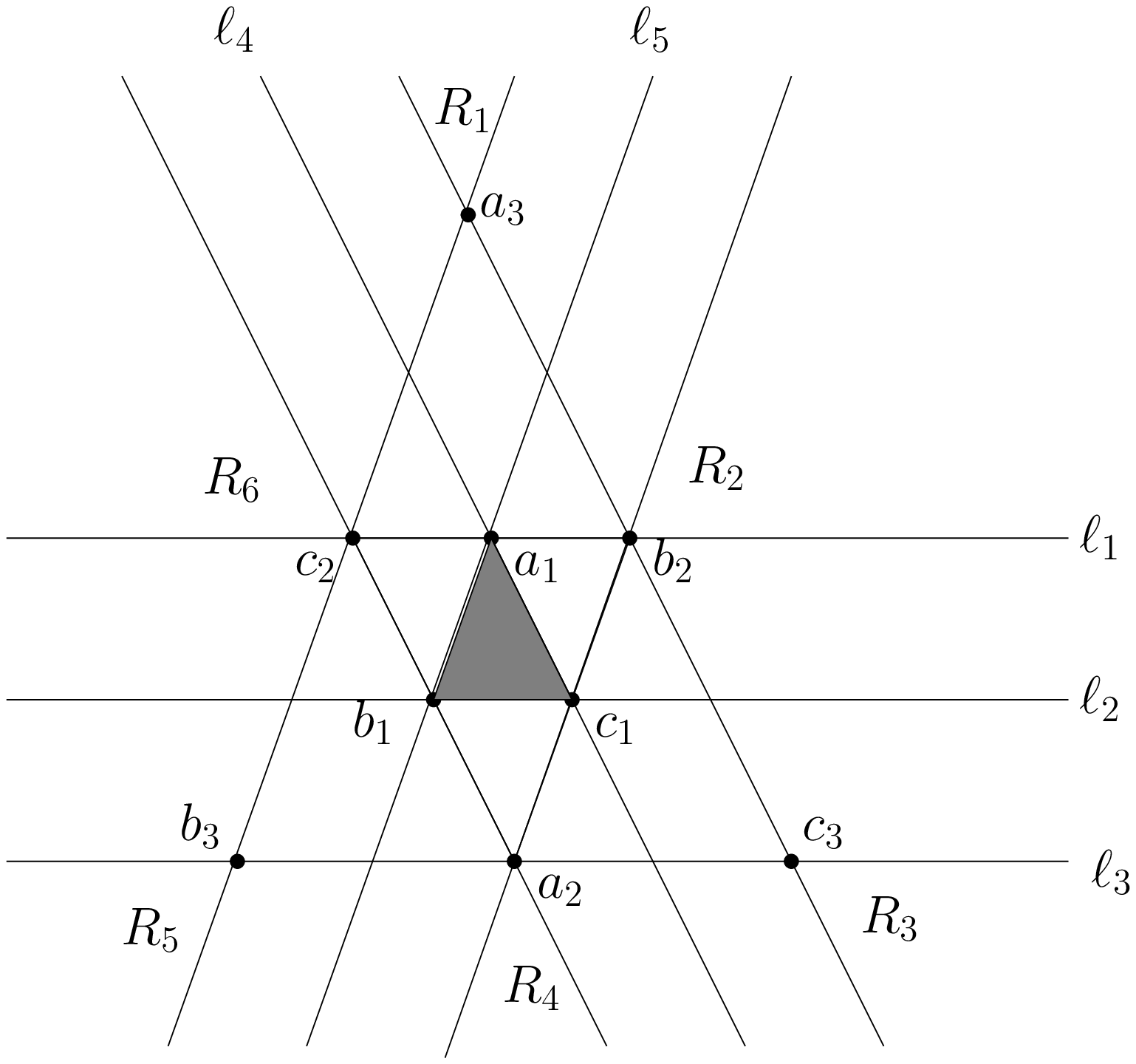}
\hskip 0.3in
\epsfxsize=3in \epsffile{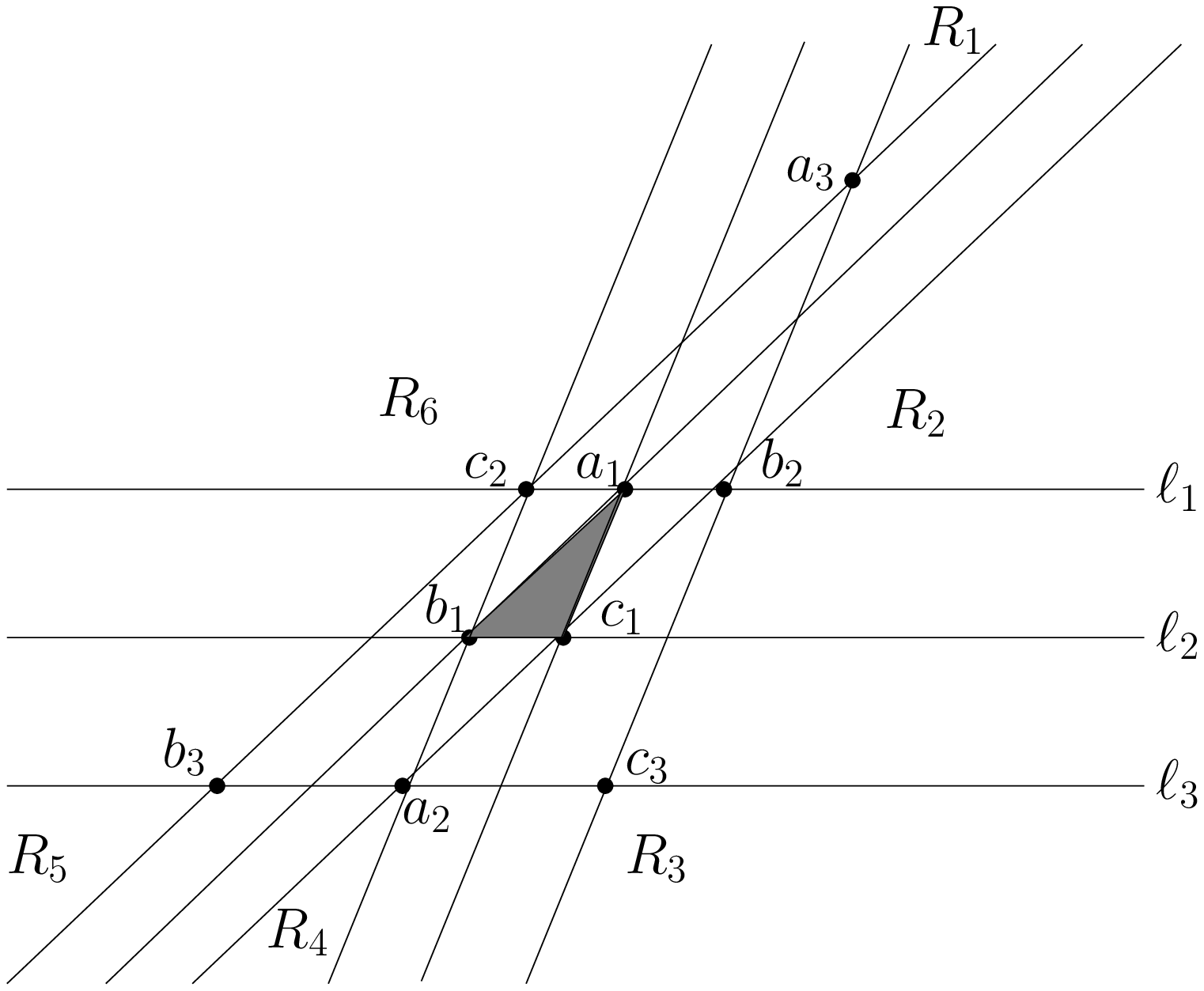}}
\centerline{\hfill (a)\hskip 2.4in (b)\hfill}
\caption{\small (a) Acute triangles: the graph $G$ is planar.
(b) Convex position: the graph $G$ is quasi-planar.}
\label{acute-obtuse}
\end{figure}

\begin{theorem} \label{thm:acute}
The maximum number of acute triangles of minimum area determined by $n$
points in the plane is $O(n)$. This bound is asymptotically tight.
\end{theorem}
\begin{proof}
Let $S$ be a set of $n$ points in the plane, and let $\T$ denote the
set of acute minimum-area triangles determined by $S$. Define a
geometric graph $G=(V,E)$ on $V=S$, where $uv \in E$ if and only if
$uv$ is a shortest side of a triangle in $\T$. We first argue that
every segment $uv$ is a shortest edge of at most two triangles in
$\T$, and then we complete the proof by showing that $G$ is planar
and so it has only $O(n)$ edges.

Let $\Delta{a_1 b_1 c_1} \in \T$ and assume that $b_1 c_1$ is a
shortest side of $\Delta{a_1 b_1 c_1}$.
Let $\Delta{a_2 b_2 c_2}$ be the triangle such that the midpoints of its
sides are $a_1,b_1,c_1$; and let $\Delta{a_3 b_3 c_3}$ be the
triangle such that the midpoints of its sides are $a_2,b_2,c_2$.
Refer to Figure~\ref{acute-obtuse}(a). Since $\Delta{a_1 b_1 c_1}$
has minimum area, then, in the notation of the figure,
each point of $S \setminus \{a_1,b_1,c_1\}$ lies in one of the
(closed) regions $R_1$ through $R_6$ or on one of the lines
$\ell_2$, $\ell_4$ or $\ell_5$; also, no point of $S \setminus
\{a_1,b_1,c_1\}$  lies in the interior of $\Delta{a_3 b_3 c_3}$.
Similarly, any point $a\in S$ of
a triangle $\Delta{ab_1c_1} \in \T$ must lie on $\ell_1$ or
$\ell_3$. Thus $a=a_1$ and $a=a_2$ are the only possible positions
of $a$. This follows from the fact that the triangles of $T$ are
acute: any point on, say, $\ell_1\cap\bd R_2$ or $\ell_1\cap\bd R_6$
forms an {\em obtuse} triangle with $b_1c_1$.

Consider two acute triangles $\Delta{a_1b_1c_1},\Delta{xyz} \in \T$
of minimum area with shortest sides $b_1 c_1\in E$ and $xy\in E$,
respectively. Assume that edges $b_1c_1$ and $xy$ cross each other.
We have the following four possibilities: (i) $x$ and $y$ lie in two
opposite regions $R_iR_{i+3}$, for some $i \in \{1,2,3\}$; (ii)
$x=a_1$ and $y \in R_4$; (iii) $x \in \ell_4$ and $y \in R_4$; (iv)
$x \in \ell_5$ and $y \in R_4$. Since $xy$ is a shortest side of
$\Delta{xyz}$, the distance from $z$ to the line through $x$ and $y$ is at
least $\sqrt{3}/2|xy|$. But then, in all four cases $\Delta{xyz}$
cannot be an acute triangle of minimum area, since it contains one
of the vertices of $\Delta{a_1 b_1 c_1}$ in its interior, a
contradiction. (For instance if $x \in R_1$ and $y \in R_4$,
$\Delta{xyc_1}$ would be obtuse and $\Delta{xyz}$ contains $c_1$ in
its interior, or if $x=a_1$ and $y \in R_4$, $\Delta{xyz}$ contains
either $b_1$ or $c_1$ in its interior.)
\end{proof}

\paragraph{Convex position.}
For points in strictly convex position we prove a tight $\Theta(n)$
bound on the maximum possible number of minimum-area triangles.
Note that a regular $n$-gon has $n$ such triangles,
so it remains to show an $O(n)$ upper bound. Also,
$n$ points equally distributed on two parallel
lines (at equal distances) give a well-known quadratic lower bound, 
so the requirement that the points be in strictly convex position
is essential for the bound to hold.

\begin{theorem} \label{thm:convex-min}
The maximum number of minimum-area triangles determined by $n$ points in
(strictly) convex position in the plane is $O(n)$. 
This bound is asymptotically tight.
\end{theorem}
\begin{proof}
The argument below is similar to that in the proof of
Theorem~\ref{thm:acute}. Since there can be only $O(n)$
acute triangles of minimum area, it is sufficient to consider
right-angled and obtuse triangles (for simplicity,
we refer to both types as obtuse), even though
the argument also works for acute triangles. We use a
similar notation: now $\T$ denotes the set of obtuse triangles of
minimum area. We define a geometric graph $G=(V,E)$ on $V=S$, where
$uv \in E$ if and only if $uv$ is a shortest side of a triangle in
$\T$. See Figure \ref{acute-obtuse}(b).

Let $\Delta{a_1 b_1 c_1} \in \T$ with $b_1c_1$ a shortest side. By
convexity, at most four triangles in $\T$ can have a common shortest
side $b_1c_1$: at most two such triangles have a third vertex on
$\ell_1$ and at most another two of them have a third vertex on $\ell_3$.
A graph drawn in the plane is said to be {\em quasi-planar}
if it has no three edges which are pairwise crossing; it is
known~\cite{aapps-97} (see also \cite{at-07})
that any quasi-planar graph with $n$ vertices
has at most $O(n)$ edges.
We now show that $G$ is quasi-planar, which will complete the proof
of the theorem.

Consider the triangles
$\Delta{a_2 b_2 c_2}$ and $\Delta{a_3 b_3 c_3}$, defined as in the
proof of Theorem~\ref{thm:acute}.
Each point of $S \setminus \{a_1,b_1,c_1\}$ lies in one of the
(closed) regions $R_1$ through $R_6$;
in particular no such point
lies in the interior of $\Delta{a_3 b_3 c_3}$. 
(Here, unlike the previous analysis, strict convexity rules out points
on any of the three middle lines, such as $\ell_2$.)
In addition, by
convexity, the regions $R_1$, $R_3$ and $R_5$ are empty of points.
Assume now that $b_1 c_1$, $xy$, $uv$ form a triplet of pairwise
crossing edges, where $xy$ and $uv$ are distinct shortest sides of
two triangles $\Delta{xyz} \in \T$ and $\Delta{uvw} \in \T$. 
It follows that each of the two edges $xy$ and $uv$ must have one
endpoint at $a_1$ and the other in $R_4$ (since each crosses $b_1 c_1$). 
Thus two edges in this triplet have  a common
endpoint, and so they do not cross, which is a contradiction.
\end{proof}

\paragraph{No three collinear points.}
We conjecture that if no three points are collinear, then the maximum
number of triangles of minimum area is close to linear. It is not
linear, though: It has been proved recently~\cite{dpt-05} that there
exist $n$-element point sets
in the plane that span $\Omega(n \log{n})$ empty congruent triangles.
Here, we show that one can repeat this construction such that there is
no collinear triples of points and that the $\Omega(n \log{n})$ empty
congruent triangles have minimum (nonzero) area.
However, we do not know of any sub-quadratic upper bound.

\begin{theorem} \label{thm:general}
  For all $n \geq 3$, there exist $n$-element point sets in the plane that
  have no three collinear points and span $\Omega(n\log n)$ triangles of
  minimum (nonzero) area.
\end{theorem}
\begin{proof}
The construction is essentially the one given in \cite{dpt-05}, and
we provide here only a brief description. We then specify the
additional modifications needed for our purposes. First, a point set
$S$ is constructed with many, i.e., $\Omega(n \log{n})$, pairwise
congruent triples of collinear points, which can be also viewed as
degenerate empty congruent triangles. Then 
this construction is slightly perturbed to obtain a set of points $S$
with no collinear triples, so that these degenerate triangles become
non-degenerate empty congruent triangles of minimum (nonzero) area. 
The details are as follows (see~\cite{dpt-05}).

Let $n=3^k$ for some $k\in \NN$.
Consider $k$ unit vectors $b_1,\ldots,b_k$, and for
$1\le i\le k$, let $\beta_i$ be the counterclockwise angle from the
$x$-axis to $b_i$.
Let $\lambda \in (0,1)$ be fixed and let $a_i=\lambda b_i$. Consider
now all $3^k$ possible sums of these $2k$ vectors, $a_i$ and $b_i$,
$1\le i\le k$,
with coefficients $0$ or $1$, satisfying the condition that for each
$i$, at least one of $a_i$ or $b_i$ has coefficient $0$.
Let $S$ be the set of $3^k$ points determined by these vectors.
Clearly, each triple of the form
($v$, $v+a_i$, $v+b_i$), where $v$ is a subset sum that does not
involve $a_i$ or $b_i$, consists of collinear points.
For such a triple, denote by $s_i(v)$ the segment whose endpoints are
$v$ and $v+b_i$.  We say that the collinear triple $(v, v+a_i, v+b_i)$
is of type $i$, $i=1,\ldots,k$.
For each $i$ there are exactly $3^{k-1}$
triples of type $i$, therefore a total of $ k 3^{k-1} = (n
\log{n})/(3 \log{3}) = \Omega(n \log{n}) $ triples of collinear
points. Clearly, all these triples form degenerate congruent
triangles in $S$.
Denote by $\ell_i(v)$ the line supporting the segment $s_i(v)$,
and by $\L$ the set of lines corresponding to these triples.

We need the following slightly stronger version of Lemma 1
in~\cite{dpt-05}. The proof is very similar to the proof of
Proposition 1 in~\cite{dpt-05}, and we omit the details.

\begin{lemma} \label{lem:exist}
There exist angles $\beta_1,\ldots,\beta_k$,
and $\lambda \in (0,1)$, such that
{\rm (i)} $S$ consists of $n$ distinct points;
{\rm (ii)} if $u,v,w \in S$ are collinear (in this order), then
$v=u+a_i$ and $w=u+b_i$.
\end{lemma}

Let $\eps$ be the minimum distance between points
$p \in S \setminus \{v, v+a_i, v+b_i\}$
and lines $\ell_i(v) \in \L$, over all pairs $(v,i)$.
By Lemma~\ref{lem:exist}, we have $\eps>0$. Now instead of choosing
$a_i$ to be collinear with $b_i$, slightly rotate $\lambda b_i$
counterclockwise from $b_i$ through a sufficiently small
angle $\delta$
about their common origin, so the collinearity disappears.
This modification is carried out at the
same time for all vectors $a_i$, $i=1,\ldots,k$,
that participate in the construction.
By continuity, there exists a sufficiently small $\delta=\delta(\eps)>0$,
so that (i) each of the triangles $\Delta(v, v+a_i, v+b_i)$
remains empty throughout this small perturbation, (ii) the point set
$S$ is in general position after the perturbation, and (iii) the
congruent triangles $\Delta(v, v+a_i, v+b_i)$ have minimum area.
This completes the proof.
\end{proof}

\vspace{-\baselineskip}

\section{Unit-area triangles in 3-space}
\label{sec:unit3}

Erd\H{o}s and Purdy~\cite{ep-71} showed that a $\sqrt{\log n}\times
(n/\sqrt{\log n})$ section of the integer lattice determines
\linebreak $\Omega(n^2 \log\log{n})$ triangles of the same area.
Clearly, this bound is also valid in 3-space. They have also
derived an upper bound of $O(n^{8/3})$ on the number of unit-area
triangles in $\RR^3$. Here we improve this bound to 
$O(n^{17/7}\beta(n)) = O(n^{2.4286})$.
We use $\beta(n)$ to denote any function of the form
$\exp(\alpha(n)^{O(1)})$, where $\alpha(n)$ is the extremely slowly
growing inverse Ackermann function. Any such function $\beta(n)$ is also
extremely slowly growing.

\begin{theorem}\label{thm:unit3}
  The number of unit-area triangles spanned by $n$ points in $\RR^3$
  is $O(n^{17/7}\beta(n)) = O(n^{2.4286})$.
\end{theorem}

The proof of the theorem is quite long, and involves several technical steps.
Let $S$ be a set of $n$ points in $\RR^3$. For each pair $a,b$ of
distinct points in $S$, let $\ell_{ab}$ denote the line passing
through $a$ and $b$, and let $C_{ab}$ denote the cylinder whose axis
is $\ell_{ab}$ and whose radius is $2/|ab|$. Clearly, any point $c\in
S$ that forms with $ab$ a unit-area triangle, must lie on $C_{ab}$.
The problem is thus to bound the number of incidences between
${n\choose 2}$ cylinders and $n$ points, but it is complicated for two
reasons: (i) The cylinders need not be distinct.  (ii) Many distinct
cylinders can share a common generator line, which may contain many
points of $S$.

\vspace{-1mm}
\paragraph{Cylinders with large multiplicity.}
Let $\C$ denote the multiset of the ${n\choose 2}$ cylinders
$C_{ab}$, for $a,b\in S$. Since the cylinders in $\C$ may appear
with multiplicity, we fix a parameter $\mu=2^j$, $j=0,1,\ldots$,
and consider separately incidences with each of the sets $\C_\mu$,
of all the cylinders whose multiplicity
is between $\mu$ and $2\mu-1$. Write $c_\mu=|\C_\mu|$.
We regard $\C_\mu$ as a set (of distinct cylinders), and will multiply
the bound that we get for the cylinders in $\C_\mu$ by $2\mu$,
to get an upper bound on the number of
incidences that we seek to estimate. We will then sum up the resulting
bounds over $\mu$ to get an overall bound.

Let $C$ be a cylinder in $\C_\mu$. Then its axis $\ell$ must contain
$\mu$ pairs of points of $P$ at a fixed distance apart (equal to $2/r$,
where $r$ is the radius of $C$). That is, $\ell$ contains
$t>\mu$ points of $S$. Let us now fix $t$ to be a power of $2$, and
consider the subset $\C_{\mu,t} \subset \C_\mu$ of those cylinders in $\C_\mu$
that have at least $t$ and at most $2t-1$ points on their axis.
By the Szemer\'edi-Trotter Theorem~\cite{st-83} (or, rather, its
obvious extension to 3-space), the number of lines containing at
least $t$ 
points of $S$ is $O(n^2/t^3 + n/t)$.
Any such line $\ell$ can be the axis of many cylinders in $\C_\mu$
(of different radii).  Any such cylinder ``charges'' $\Theta(\mu)$
pairs of points out of the $O(t^2)$ pairs along $\ell$,
and no pair is charged more than once.
Hence, for a given line $\ell$ incident to at least $t>\mu$ and at
most $2t-1$ points of $S$, the number of distinct cylinders in
$\C_\mu$ that have $\ell$ as axis is $O(t^2/\mu)$. Summing over
all axes incident to at least $t$ and at most $2t-1$ points
yields that the number of distinct cylinders in $\C_{\mu,t}$ is
\begin{equation}\label{nut}
c_{\mu,t}=
O\left(\left(\frac{n^2}{t^3} + \frac{n}{t}\right)\frac{t^2}{\mu}\right)=
O\left(\frac{n^2}{t\mu} + \frac{nt}{\mu}\right).
\end{equation}
We next sum this over $t$, a power of 2 between $\mu$ and $\nu$,
and conclude that the number of distinct cylinders
in $\C_\mu$ having at most $\nu$ points on their axis is
\begin{equation} \label{cmu}
c_{\mu,\le\nu} = O\left( \frac{n^2}{\mu^2} + \frac{n \nu}{\mu} \right).
\end{equation}

\vspace{-1mm}
\paragraph{Restricted incidences between points and cylinders.}
We distinguish two {\em type}s of incidences, which we count
separately. An incidence between a point $p$ and a cylinder $C$ is of
{\em type 1} if the generator of $C$ passing through $p$ contains at
least one additional point of $S$; otherwise it is of {\em type 2}. We
begin with the following subproblem, in which we bound the number of
incidences between the cylinders of $\C$, counted with multiplicity,
and {\em multiple} points that lie on their generator lines, as well
as incidences with cylinders with ``rich'' axes.
Specifically, we have the following lemma.
\begin{lemma} \label{107}
  Let $S$ be a set of $n$ points and $\C$ be the multiset of the
  ${n\choose 2}$ cylinders $C_{ab}$, for $a,b\in S$ (counted with
  multiplicity). The
  total number of all incidences of type 1 and all incidences
  involving cylinders having at least $n^{14/45}$ points on their axis
  is bounded by $O(n^{107/45}\polylog(n))=O(n^{2.378})$.
\end{lemma}
\begin{proof}
Let $L$ denote the set of lines spanned by the points of $S$. Fix
a parameter $k=2^i$, $i=1,\ldots$,
and consider the set $L_k$ of all lines that contain at least $k$ and
at most $2k-1$ points of $S$.
We bound the number of incidences between
cylinders in $\C$ that contain lines in $L_k$ as generators
and points that lie on those lines. Formally, we bound the
number of triples $(p,\ell,C)$, where $p\in S$, $\ell\in L_k$,
and $C\in\C$, such that $p\in\ell$ and $\ell\subset C$.
Summing these bounds over $k$ will give us a bound for the
number of incidences of type 1. Along the way, we will also dispose of
incidences with cylinders whose axes contain many points.

As already noted, the Szemer\'edi-Trotter Theorem~\cite{st-83} implies that
${\displaystyle
\lambda_k:=|L_k|=O\left( \frac{n^2}{k^3} + \frac{n}{k} \right)}$.

\vspace{-1mm}
\paragraph{Line-cylinder incidences.}
Consider the subproblem of bounding the number of
incidences between lines in $L_k$ and cylinders in $\C$, where
a line $\ell$ is said to be incident to cylinder $C$ if $\ell$
is a generator of $C$. We will then multiply the resulting
bound by $2k$ to get an upper bound on the number of
point-line-cylinder incidences involving $L_k$, and then
sum the resulting bounds over $k$.

\vspace{-1mm}
\paragraph{Generator lines with many points.}
Let us first dispose of the case $k>n^{1/3}$. Any line
$\ell\in L_k$ can be a generator of at most $n$ cylinders (counted
with multiplicity), because, having fixed $a\in S$, the point
$b\in S$ such that $C_{ab}$ contains $\ell$ is determined
(up to multiplicity $2$). Hence the number of incidences between
the points that lie on $\ell$ and the cylinders of $\C$ is
$O(nk)$. Summing over $k=2^i>n^{1/3}$ yields the overall bound
$$
O\left( \sum_k nk\lambda_k \right) =
O\left( \sum_k \left( \frac{n^3}{k^2} + n^2 \right) \right) = O(n^{7/3}) .
$$
Hence, in what follows, we may assume that $k\le n^{1/3}$.
In this range of $k$ we have
\begin{equation} \label{lambdak}
\lambda_k=O\left( \frac{n^2}{k^3} \right) .
\end{equation}

\vspace{-1mm}
\paragraph{Axes with many points.}
Let us also fix the multiplicity $\mu$ of the cylinders under
consideration (up to a factor of $2$, as above).
The number of distinct cylinders in $\C_\mu$
having between $t>\mu$ and $2t-1$ points on their axes, is
$O(n^2/(t\mu) + nt/\mu)$; see (\ref{nut}). While the first term
is sufficiently small for our purpose, the second
term may be too large when $t$ is large. To avoid this difficulty, we
fix another threshold exponent $z < 1/2$ that we will optimize later,
and handle separately the cases $t\ge n^{z}$ and $t< n^{z}$.
That is, in the first case, for $t\ge n^{z}$ a power of 2, we seek an
upper bound on the overall number of incidences between the points of
$S$ and the cylinders in $\C$ whose axis contains between $t$ and $2t-1$
points of $S$. (For this case, we combine all the multiplicities
$\mu<t$ together.) By the Szemer\'edi-Trotter theorem~\cite{st-83},
the number of such axes is $O(n^2/t^3+n/t)$.

Fix such an axis $\alpha$. It defines $\Theta(t^2)$ cylinders,
and the multiplicity of any of these cylinders is at most $O(t)$. Since no
two distinct cylinders in this collection can pass through the same
point of $S$, it follows that the total number of incidences
between the points of $S$ and these cylinders is $O(nt)$.
Hence the overall number of incidences under consideration is
$O(n^2/t^3+n/t) \cdot O(nt)= O(n^3/t^2+n^2)$.
Summing over all $t\ge n^z$,
a power of $2$, we get the overall bound $O(n^{3-2z})$.

Note that this bound takes care of {\em all} the incidences between
the points of $S$ and the cylinders having at least $t\ge n^z$ points
along their axes, not just those of type 1 (involving multiple points
on generator lines).

\paragraph{Cylinders with low multiplicity.}
We now confine the analysis to cylinders having fewer than $n^z$
points on their axis, and go back to fixing the multiplicity $\mu$,
which we may assume to be at most $n^z$.
We thus want to bound the number of
incidences between $\lambda_k$ distinct lines and $c_{\mu,\le n^z}$
distinct cylinders in $\C_\mu$, for given $k\le n^{1/3}$, $\mu \leq n^z$.
Note that a cylinder can contain a line if and only if it is parallel to the
axis of the cylinder, so we can split the
problem into subproblems, each associated with some direction
$\theta$, so that in the $\theta$-subproblem we have a set of some
$c_\mu^{(\theta)}$ cylinders and a set of some
$\lambda_k^{(\theta)}$ lines, so that the lines and the cylinder axes
are all parallel (and have direction $\theta$);
we have $\sum_\theta c_\mu^{(\theta)} = c_{\mu,\le n^z}$,
and $\sum_\theta \lambda_k^{(\theta)} = \lambda_k$.

For a fixed $\theta$, we project the cylinders and lines in the
$\theta$-subproblem onto a plane with normal direction $\theta$,
and obtain a set of $c_\mu^{(\theta)}$ circles and a set of
$\lambda_k^{(\theta)}$ points, so that the number of line-cylinder
incidences is equal to the number of point-circle incidences.
By \cite{ANPPSS,ArS,MT},\footnote{%
  The bound that we use, from \cite{MT}, is slightly better than
  the previous ones.}
the number of point-circle incidences between $N$ points and $M$ circles in the plane is
$O(N^{2/3}M^{2/3} + N^{6/11}M^{9/11}\log^{2/11}(N^3/M)+N+M)$.
It follows that the number of such line-cylinder incidences is
\begin{equation} \label{eq:theta}
O\left(
(\lambda_k^{(\theta)})^{2/3}
(c_\mu^{(\theta)})^{2/3} +
(\lambda_k^{(\theta)})^{6/11}
(c_\mu^{(\theta)})^{9/11}
  \log^{2/11}((\lambda_k^{(\theta)})^3/c_\mu^{(\theta)})+
\lambda_k^{(\theta)} +
c_\mu^{(\theta)} \right) .
\end{equation}
Note that, for any fixed $\theta$, we have
$\lambda_k^{(\theta)} \le n/k$ and
$c_\mu^{(\theta)} \le n^{1+z}/\mu$.
The former inequality is trivial.  To see the latter inequality,
note that an axis with $t<n^z$ points
defines ${t\choose 2}$ cylinders. Since we only consider cylinders
with multiplicity $\Theta(\mu)$, the number of distinct such
cylinders is $O(t^2/\mu)$, and the number of lines
(of direction $\theta$) with about $t$
points on them is at most $n/t$, for a total of at most $O(nt/\mu)$
distinct cylinders. Partitioning the range $\mu < t\le n^{z}$ by
powers of $2$, as above, and summing up the resulting bounds, the bound
$c_\mu^{(\theta)} \le n^{1+z}/\mu$ follows.

Summing over $\theta$, and using H\"older's inequality, we have
(here $x$ is a parameter between $2/11$ and $6/11$
that we will fix shortly)
$$
\sum_\theta
(\lambda_k^{(\theta)})^{6/11}
(c_\mu^{(\theta)})^{9/11} \le
\left(\frac{n}{k}\right)^{6/11-x}
\left(\frac{n^{1+z}}{\mu}\right)^{x-2/11}
\sum_\theta
(\lambda_k^{(\theta)})^{x}
(c_\mu^{(\theta)})^{1-x} \le
$$
$$
\frac{n^{(4-2z)/11+xz}}{k^{6/11-x}\mu^{x-2/11}}
\left( \sum_\theta
\lambda_k^{(\theta)} \right)^{x}
\left( \sum_\theta
c_\mu^{(\theta)} \right)^{1-x} =
\frac{n^{(4-2z)/11+xz}}{k^{6/11-x}\mu^{x-2/11}}
\lambda_k^x c_{\mu,\le n^z}^{1-x} .
$$
We need to multiply this bound by $\Theta(k\mu)$.
Substituting the bounds $\lambda_k=O(n^2/k^3)$ from (\ref{lambdak}),
and $c_{\mu,\le n^z} = O(n^2/\mu^2 + n^{1+z}/\mu)$ from (\ref{cmu}),
we get the bound
\begin{eqnarray*}
&& O\left(
n^{(4-2z)/11+xz} k^{5/11+x}\mu^{13/11-x}
\left( \frac{n^2}{k^3} \right)^x
\left( \frac{n^2}{\mu^2} + \frac{n^{1+z}}{\mu} \right)^{1-x}
\log^{2/11} n\right)\\
&=&O\left( k^{5/11-2x} \left(
n^{2+(4-2z)/11+xz} \mu^{x-9/11} +
n^{(15+9z)/11+x} \mu^{2/11} \right)
\log^{2/11} n \right) .
\end{eqnarray*}
Choosing $x=5/22$ (the smallest value for which the exponent of $k$
is non-positive), the first term becomes
$O( n^{2+4/11+z/22} \log^{2/11} n)$, which we need to balance with
$O(n^{3-2z})$; for this, we choose $z=14/45$ and obtain the bound
$O(n^{107/45}\log^{2/11}n)=O(n^{2.378})$; for this choice of $z$,
recalling that $\mu<n^z$, the second term is dominated by the first.
Summing over $k,\mu$ only adds logarithmic factors, for a resulting
overall bound $O(n^{2.378})$.

Similarly, we have (with a different choice of $x$, soon to be made)
$$
\sum_\theta
(\lambda_k^{(\theta)})^{2/3}
(c_\mu^{(\theta)})^{2/3} \le
\left(\frac{n}{k}\right)^{2/3-x}
\left(\frac{n^{1+z}}{\mu}\right)^{x-1/3}
\sum_\theta
(\lambda_k^{(\theta)})^{x}
(c_\mu^{(\theta)})^{1-x} \le
$$
$$
\frac{n^{(1-z)/3+xz}}{k^{2/3-x}\mu^{x-1/3}}
\left( \sum_\theta
\lambda_k^{(\theta)} \right)^{x}
\left( \sum_\theta
c_\mu^{(\theta)} \right)^{1-x} =
\frac{n^{(1-z)/3+xz}}{k^{2/3-x}\mu^{x-1/3}}
\lambda_k^x c_{\mu,\le n^z}^{1-x} .
$$
Multiplying by $k\mu$ and arguing as above, we get
\begin{eqnarray*}
&&O\left(
n^{(1-z)/3+xz} k^{1/3+x}\mu^{4/3-x}
\left( \frac{n^2}{k^3} \right)^x
\left( \frac{n^2}{\mu^2} + \frac{n^{1+z}}{\mu} \right)^{1-x}
\log^{2/11} n\right)\\
&=&
O\left( k^{1/3-2x} \left(
n^{2+(1-z)/3+xz} \mu^{x-2/3} +
n^{1+(1+2z)/3+x} \mu^{1/3} \right)
\log^{2/11} n\right) .
\end{eqnarray*}
We choose here $x=1/6$ and note that, for $z=14/45$ and $\mu < n^z$,
the bound is smaller than $O(n^{7/3})$, which is dominated by the
preceding bound $O(n^{2.378})$.

Finally, the linear terms in (\ref{eq:theta}), multiplied by $k\mu$,
add up to 
$$
k\mu \sum_{\theta} O\left(\lambda_k^{(\theta)}+c_\mu^{(\theta)}\right)=
O\left( k\mu \left( \lambda_k + c_{\mu,\le n^z} \right)\right) =
O\left( \frac{n^2\mu}{k^2} + \frac{n^2k}{\mu} + n^{1+z}k \right) ,
$$
which, by our assumptions on $k$, $\mu$, and $z$ is also dominated by
$O(n^{2.378})$. Summing over $k,\mu$ only add logarithmic factors, for
a resulting overall bound $O(n^{2.378})$.
This completes the proof of Lemma~\ref{107}.
\end{proof}

It therefore remains to count point-cylinder incidences of type 2,
involving cylinders having at most $n^{14/45}$ points on their axes.

\vspace{-1mm}
\paragraph{The intersection pattern of three cylinders.}
We need the following technical lemma, whose proof is borrowed from a
yet unpublished work \cite{EPS:cyl}, and is presented in the appendix.

\begin{lemma}\label{3int}
  Let $C,C_1,C_2$ be three cylinders with no pair of parallel axes.
  Then $C\cap C_1\cap C_2$ consists of at most 8 points.
\end{lemma}

\paragraph{Point-cylinder incidences.}
Using the partition technique~\cite{ch-05,ps-04}
for disjoint cylinders in $\RR^3$, we show the following:

\begin{lemma}\label{lem:cut}
  For any parameter $r$, $1\leq r\leq \min\{m,n^{1/3}\}$,
  the maximum number of incidences
  of type 2 between $n$ points and $m$ cylinders in 3-space satisfies
  the following recurrence:
\begin{equation}\label{eq:inci}
  I(n,m)=O(n+mr^2\beta(r)) + O(r^3\beta(r))\cdot
  I\left(\frac{n}{r^3},\frac{m}{r}\right),
\end{equation}
for some slowly growing function $\beta(n)$, as above.
\end{lemma}
\begin{proof}
Let $\C$ be a set of $m$ cylinders, and $S$ be a set of $n$ points.
Construct a $(1/r)$-cutting of the arrangement $\A(\C)$.  The
cutting has $O(r^3\beta(r))$ relatively open pairwise disjoint
cells, each crossed by at most $m/r$ cylinders
and containing at most $n/r^3$ points of $S$~\cite{ceg-89}
(see also~\cite[p.~271]{as-95}); the first property is by definition
of $(1/r)$-cuttings, and the second is enforced by subdividing cells
with too many points.  The number of incidences between
points and cylinders {\em crossing} their cells is thus
$$
  O(r^3\beta(r))\cdot I\left(\frac{n}{r^3},\frac{m}{r}\right) .
$$
(Note that any incidence of type 2 remains an incidence of type 2 in
the subproblem it is passed to.)

It remains to bound the number of incidences between the points of
$S$ and the cylinders that {\em contain} their cells. Let $\tau$ be
a (relatively open) lower-dimensional cell of the cutting.
If ${\rm dim}(\tau)=2$ then we can assign any point $p$ in $\tau$
to one of the two neighboring full-dimensional cells, and count
all but at most one of the incidences with $p$ within that cell.
Hence, this increases the count by at most $n$.

If ${\rm dim}(\tau)=0$, i.e., $\tau$ is a vertex of the cutting,
then any cylinder containing $\tau$ must cross or define one of the
full-dimensional cells adjacent to $\tau$. Since each cell has at
most $O(1)$ vertices, it follows that the total number of such
incidences is $O(r^3\beta(r))\cdot (m/r) = O(mr^2\beta(r))$.

Suppose then that ${\rm dim}(\tau)=1$, i.e., $\tau$ is an edge of
the cutting. An immediate implication of Lemma~\ref{3int} is that
only $O(1)$ cylinders can contain $\tau$, unless $\tau$ is a line,
which can then be a generator of arbitrarily many cylinders.

Since we are only counting incidences of type 2, this implies that
any straight-edge 1-dimensional cell $\tau$ of the cutting generates
at most one such incidence with any cylinder that fully contains
$\tau$.  Non-straight edges of the cutting are contained in only
$O(1)$ cylinders, as just argued, and thus the points on such edges
generate a total of only $O(n)$ incidences with the cylinders. 
Thus the overall number of incidences in
this subcase is only $O(n+r^3\beta(r))$. Since $r\le m$,
this completes the proof of the lemma.
\end{proof}

\begin{lemma}\label{lem:distinct-cylinders}
  The number of incidences of type 2 between $n$ points and $m$
  cylinders in $\RR^3$ is
\begin{equation} \label{177}
O\left(\left( m^{6/7}n^{5/7}+m+n \right) \beta(n)\right).
\end{equation}
\end{lemma}
\begin{proof}
Let $\C$ be a set of $m$ cylinders, and $S$ be a set of $n$ points.
We first derive an upper bound of $O(n^5+m)$
on the number of incidences of type 2 between $\C$ and $S$. We represent
the cylinders as points in a dual 5-space, so that each cylinder
$C$ is mapped to a point $C^*$, whose coordinates are the five
degrees of freedom of $C$ (four specifying its axis and the fifth
specifying its radius). A point $q\in \reals^3$ is mapped to a
surface $q^*$ in $\reals^5$, which is the locus of all points dual
to cylinders that are incident to $q$. With an appropriate choice of
parameters, each surface $q^*$ is semi-algebraic of constant description
complexity. By definition, this duality preserves incidences.

After dualization, we have an incidence problem involving $m$ points
and $n$ surfaces in $\reals^5$.  We construct the arrangement
$\A$ of the $n$ dual surfaces, and bound the number of their
incidences with the $m$ dual points as follows.  The arrangement
$\A$ consists of $O(n^5)$ relatively open cells of dimensions
$0,1,\ldots,5$. Let $\tau$ be a cell of $\A$. We may assume that
${\rm dim}(\tau)\le 4$, because no point in a full-dimensional cell
can be incident to any surface.

If $\tau$ is a vertex, consider any surface $\varphi$ that passes
through $\tau$. Then $\tau$ is a vertex of the arrangement
restricted to $\varphi$, which is a $4$-dimensional arrangement with
$O(n^4)$ vertices. This implies that the number of incidences at
vertices of $\A$ is at most $n\cdot O(n^4) = O(n^5)$.

Let then $\tau$ be a cell of $\A$ of dimension $\ge 1$, and let $u$
denote the number of surfaces that contain $\tau$. If $u\le 8$ then each
point in $\tau$ (dual to a cylinder) has at most $O(1)$ incidences
of this kind, for a total of $O(m)$.

Otherwise, $u\ge 9$. Since ${\rm dim}(\tau)\ge 1$, it contains
infinitely many points dual to cylinders (not necessarily in $\C$).
By Lemma~\ref{3int}, back in the primal 3-space, if three cylinders
contain the same nine points, then the axes of at least two of
them are parallel. Hence all $u$ points lie on one line or on two
parallel lines, which are common generators of these pair of
cylinders. In this case, all cylinders whose dual points lie in $\tau$
contain these generator(s). But then, by definition, the incidences
between these points and the cylinders of $\C$ whose dual points lie
on $\tau$ are of type 1, and are therefore not counted at all by
the current analysis. Since $\tau$ is a face of $\A$, no other point
lies on any of these cylinders, so we may ignore them completely.

Hence, the overall number of incidences under consideration is $O(n^5+m)$.

If $m>n^5$, this bound is $O(m)$. If $m<n^{1/3}$,
we apply Lemma \ref{lem:cut} with $r=m$, which then yields that
each recursive subproblem has at most one cylinder, so each point in a
subproblem generates at most one incidence, for a total of $O(n)$
incidences. Hence, in this case (\ref{eq:inci}) implies that the
number of incidences between $\C$ and $S$ is
$O(n+m^3\beta(m)) = O(n\beta(n))$.

Otherwise we have $n^{1/3} \leq m \leq n^5$, so we can apply
Lemma~\ref{lem:cut} with parameter $r=(n^5/m)^{1/14}$; observe that
$1 \leq r \leq \min\{m, n^{1/3}\}$ in this case. Using the above bound
for each of the subproblems in the recurrence, we obtain
$I(n/r^3,m/r)=O((n/r^3)^5+m/r)$, and thus the total number of
incidences of type 2 in this case is
$$
  O(n+mr^2\beta(r)) + O(r^3\beta(r))\cdot O\left(
    \left(\frac{n}{r^3}\right)^5 + \frac{m}{r}\right) = O\left(
    \frac{n^5}{r^{12}} + mr^2 \right)\beta(r).
$$
The choice $r=(n^5/m)^{1/14}$ yields the bound (\ref{177}).
Combining this with the other cases, the bound in the lemma follows.
\end{proof}

We are now in position to complete the proof of
Theorem~\ref{thm:unit3}.

\begin{proofof}{Theorem~\ref{thm:unit3}}
We now return to our original setup, where the cylinders in $\C$ may
have multiplicities. We fix some parameter $\mu$ and consider,
as above, all cylinders in $C_\mu$, and recall our choice of $z=14/45$.
The case $\mu\ge n^z$ is taken care of by Lemma~\ref{107}, accounting
for at most  $O(n^{107/45}\polylog(n))$ incidences.
In fact, Lemma~\ref{107} takes
care of all cylinders that contain at least $n^z$ points on
their axes. Assume then that $\mu < n^z$, and consider only those cylinders
in $\C_\mu$ containing fewer than $n^z$ points on their axes.
By (\ref{cmu}), we have $c_{\mu,\le n^z} = O(n^2/\mu^2)$.
Consequently, the number of incidences with the remaining cylinders
in $C_\mu$, counted with multiplicity, but excluding multiple points
on the same generator line, is
$$
  O\left(\mu\beta(n)\cdot \left( \left(\frac{n^2}{\mu^2}\right)^{6/7}
   \cdot   n^{5/7} + \frac{n^2}{\mu^2} + n \right) \right)
  = O\left(\left(\frac{n^{17/7}}{\mu^{5/7}}
    + \frac{n^2}{\mu} + n\mu \right) \beta(n)\right) .
$$
Summing over all $\mu\leq n^z$ (powers of 2), and adding the
bound $O(n^{107/45}\polylog(n))=O(n^{2.378})$ from Lemma~\ref{107} on
the other kinds of incidences, we get the desired overall bound of
$O(n^{17/7}\beta(n))=O(n^{2.4286})$.
\end{proofof}

\noindent
{\bf Remark.}
In a nutshell, the ``bottleneck'' in the analysis is the case where
$\mu$ is small (say, a constant) and we count incidences of type 2.
The rest of the analysis, involved as it is, just shows that all the
other cases contribute fewer (in fact, much fewer) incidences. One
could probably simplify some parts of the analysis, at the cost of
weakening the other bounds, but we leave these parts as they are, in
the hope that the bottleneck case could be improved, in which case
these bounds might become the dominant ones.

\section{Minimum-area triangles in 3-space}
\label{sec:minimum3}

Place $n$ equally spaced
points on the three parallel edges of a right prism whose base is an
equilateral triangle, such that inter-point distances are small along
each edge. This construction yields $\frac{2}{3}n^2-O(n)$ minimum-area triangles,
a slight improvement over the lower bound construction in the plane.
Here is yet another construction with the same constant $2/3$ in the leading
term: Form a rhombus in the $xy$-plane from two equilateral triangles
with a common side, extend it to a prism in 3-space, and place $n/3$
equally spaced points on each of the lines passing through the vertices
of the shorter diagonal of the rhombus, and $n/6$ equally spaced
points on each of the two other lines, where again the
inter-point distances along these lines are all equal and small.
The number of minimum-area triangles is
$$ 2\left(\frac{1}{3 \cdot 3} + \frac{4}{3 \cdot 6}\right)n^2 -O(n)=
\frac{2}{3}n^2-O(n). $$
The following theorem shows that this bound is optimal up to a
constant factor. No quadratic upper bound has previously been known for
minimum-area triangles in $\RR^3$.

\begin{theorem} \label{thm:min3}
  The number of triangles of minimum (nonzero) area spanned by $n$ points in
  $\RR^3$ is at most $n^2+O(n)$.
\end{theorem}
\begin{proof}
Consider a set $S$ of $n$ points in $\RR^3$, and let $T$ be the set of
triangles of minimum (nonzero) area spanned by $S$.
Without loss of generality, assume the minimum area to be $1$.
Similarly to the planar case, we assign each
triangle in $T$ to one of its longest sides, and argue that at most a
constant number of triangles are assigned to each segment spanned by
$S$. This immediately implies an upper bound of $O(n^2)$ on the
cardinality of $T$. To improve the main coefficient in this bound, we
distinguish between {\em fat} and {\em thin} triangles.
A triangle is called fat (resp., thin) if the length of the height
corresponding to its longest side is at least (resp., less than)
half of the length of the longest side.
We show that the number $N_1$ of thin triangles of minimum area is
at most $2{n \choose 2}=n^2-n$, and that
the number $N_2$ of fat triangles of minimum area is only $O(n)$.

Consider a segment $ab$, with $a,b\in S$, and let $h=|ab|$.  Every
point $c\in S\setminus \{a,b\}$ for which the triangle
$\Delta{abc}$ has minimum (unit) area must lie on a bounded cylinder $C$
with axis $ab$, radius $r=2/h$,
and bases that lie on the planes $\pi_a$ and $\pi_b$, incident to
$a$ and $b$, respectively, and orthogonal to $ab$.  In fact, if
$\Delta abc$ is assigned to $ab$ (that is, $ab$ is the longest side), 
then $c$ must lie on a smaller
portion $C'$ of $C$, bounded by bases that intersect $ab$ at points at
distance $h-\sqrt{h^2-r^2}$ from $a$ and $b$, respectively.
Assume for
convenience that $ab$ is vertical, $a$ is the origin and
$b=(0,0,h)$. Since $ab$ is the longest side of $\Delta abc$, the
side of the isosceles triangle with base $ab$ and height $r$ must
be no larger than $h$, i.e., $\frac14 h^2+r^2 \le h^2$, or
$r^2\le \frac34 h^2$.
Notice that the triangle formed by any two points
of $S$ lying on $C'$ with either $a$ or $b$ is non-degenerate.

We first derive a simple formula that relates the area of any
(slanted) triangle to the area of its $xy$-projection. Consider a
triangle $\Delta$ that is spanned by two vectors $u,v$, and let
$\Delta_0$, $u_0$, and $v_0$ denote the $xy$-projections of $\Delta$,
$u$, and $v$, respectively. Write (where $\kk$ denotes, as usual, the
vector $(0,0,1)$)
$$
u=u_0 + x\kk \quad\text{and}\quad v=v_0 + y\kk ,
$$
and put $A= {\rm area}(\Delta)$, $A_0= {\rm area}(\Delta_0)$.
Then
$$
A^2 = \frac14 \|u\times v\|^2 =
\frac14\| (u_0 + x\kk)\times (v_0 + y\kk) \| =
\frac14\left( \| u_0\times v_0 \|^2 + \|yu_0 - xv_0\|^2 \right) ,
$$
or
\begin{equation} \label{eq:pa}
A^2 = A_0^2 + \frac14 \|yu_0 - xv_0\|^2 .
\end{equation}

\vspace{-1mm}
\paragraph{An initial weaker bound.}
We claim that at most 10 triangles are assigned to $ab$. Assume, to
the contrary, that this number is at least 11.  Divide $C$ into two
equal slices by a horizontal plane orthogonal to $ab$ through its
midpoint. Since more than 10 points of $S$ lie on $C$, at least
6 of them must lie on the same slice $C_0$, say the bottom slice.
It follows that two points, $c$ and $d$,
lie in some sector $\Upsilon$ of $C_0$ making a dihedral angle
$\alpha$ at $ab$ of at most $360^\circ/6=60^\circ$.
An illustration is provided in Figure~\ref{fig:min3}.

\begin{figure}[htbp]
  \begin{center}
\epsfig{file=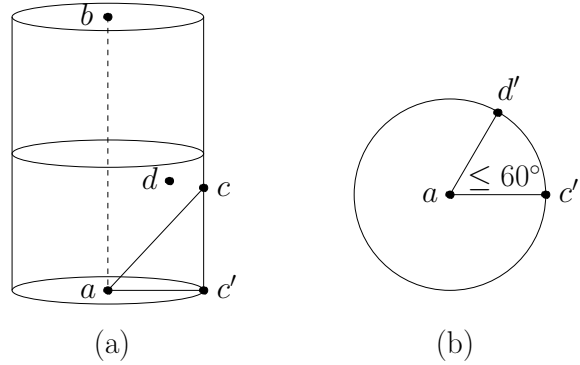,width=3in,clip=}
 \caption{\small Charging scheme for minimum-area triangles in
   3-space; (a) the cylinder $C$; (b) the projection on $\pi_a$;
   $c'$ and $d'$ are the respective projections of $c$ and $d$.}
\end{center}
\label{fig:min3}
\end{figure}

We may assume, without loss of generality, that
$$
c=(r,0,x)=c_0+x\kk \quad\text{and}\quad
d=(r\cos\alpha,r\sin\alpha,y)=d_0+y\kk ,
$$
where $0\le\alpha\le 60^\circ$ and $0\le x,y\le h/2$.
Write $A={\rm area}(\Delta{acd})$. Using (\ref{eq:pa}), we have
$$
A^2=
\frac14 \left\| c_0\times d_0\right\|^2 + \frac14\left\|yc_0-xd_0\right\|^2 =
\frac{r^4\sin^2\alpha}{4} + \frac{r^2}{4}\left( x^2+y^2-2xy\cos\alpha\right).
$$
The expression $x^2+y^2-2xy\cos\alpha$ is the squared length of the
third side of the triangle with sides $x$, $y$, with the angle
$\alpha \leq 60^\circ$ between them. Since $x,y \leq h/2$, we clearly
have $x^2+y^2-2xy\cos\alpha \leq h^2/4$.
Thus, recalling that $r^2\le \frac34h^2$ and that $h^2r^2 = 4$, we have
$$
A^2 \le \frac{r^4\sin^2\alpha}{4} + \frac{r^2}{4} \cdot \frac{h^2}{4}=
\frac{r^2}{4} \left( r^2\sin^2\alpha+ \frac{h^2}{4} \right)
\leq \frac{r^2 h^2}{4} \left( \frac{9}{16} + \frac{1}{4} \right)=
\frac{13}{16} <1,
$$
which contradicts the minimality of the area of $\Delta abc$.
Hence, at most 10 triangles are assigned to each segment
spanned by $S$.
This already implies that there are at most $5(n^2-n)$
minimum-area triangles.

\vspace{-1mm}
\paragraph{A better bound.}
We now improve the constant of proportionality, using a more careful
analysis, which distinguishes between the cases in which the
minimum-area triangles charged to the segment $ab$ are thin or fat.

\smallskip
\noindent{\bf (a)} $r< \frac12h$ (thin triangles).
We claim that in this case at most two triangles can be assigned
to $ab$. Indeed, suppose to the contrary that at least three triangles
are assigned to $ab$, so their third vertices, $c,d,e \in S$ lie on
$C'\subset C$. 
Write the $z$-coordinates of $c,d,e$ as $z_1h$, $z_2h$, $z_3h$,
respectively, and assume, without loss of generality, that
$0< z_1\le z_2\le z_3< 1$, and $z_2\le 1/2$.  Consider the
triangle $\Delta{acd}$, and let $A$ denote its area.
As before, write, without loss of generality,
$$
c=(r,0,z_1h) \quad\text{and}\quad
d=(r\cos\alpha,r\sin\alpha,z_2h) ,
$$
for some $0\le\alpha\le 180^\circ$.  Using (\ref{eq:pa}), we get
$$
A^2=\frac14 r^4\sin^2\alpha + \frac14 r^2h^2 (z_1^2+z_2^2-2z_1z_2\cos\alpha).
$$
Thus, recalling that $r< \frac12h$ and that $h^2r^2 = 4$, we get
\begin{equation} \label{eq:A2}
A^2 <
\frac14 r^2h^2 \left( \frac14\sin^2\alpha +
z_1^2+z_2^2-2z_1z_2\cos\alpha \right) =
\frac14\sin^2\alpha + z_1^2+z_2^2-2z_1z_2\cos\alpha .
\end{equation}
Let us fix $z_1,z_2$ and vary only $\alpha$. Write
$$
f(\alpha) =
\frac14\sin^2\alpha + z_1^2+z_2^2 - 2z_1z_2\cos\alpha ,
\quad\text{and}\quad
f'(\alpha) =
\frac12\sin\alpha\cos\alpha + 2z_1z_2\sin\alpha .
$$
$f$ attains its maximum at the zero of its derivative, namely at
$\alpha_0$ that satisfies
$$
\cos\alpha_0 = - 4z_1z_2 .
$$
(Note that since $z_1\le z_2\le\frac12$, we always have $4z_1z_2\le 1$.
Also, at the other zero $\alpha=0$, $f$ attains its minimum
$(z_1-z_2)^2$.)

Substituting $\alpha_0$ into (\ref{eq:A2}),
and using $z_1\le z_2\le\frac12$, we get
$$
A^2 < \frac{1-16z_1^2z_2^2}{4} + z_1^2 + z_2^2 + 8z_1^2z_2^2 =
\frac14 + z_1^2 + z_2^2 + 4z_1^2z_2^2 =
\left(\frac12 + 2z_1^2\right) \left(\frac12 + 2z_2^2\right) \le 1 ,
$$
which contradicts the minimality of the area of $\Delta abc$
(recall that $\Delta acd$ is non-degenerate).

We have thus shown that at most two thin triangles of minimum area can be
assigned to any segment $ab$, so $N_1 \leq 2{n \choose 2}=n^2-n$.

\smallskip
\noindent{\bf (b)} $r\ge \frac12h$ (fat triangles).
Recall that we always have $r\le \frac{\sqrt{3}}{2}h$. Multiplying
these two inequalities by $h/2$, we get
$$
\frac{h^2}{4} \le 1 \le \frac{h^2\sqrt{3}}{4} ,
\quad\text{or}\quad
\frac{2}{3^{1/4}} \le h \le 2 .
$$
Let $E$ denote the set of all segments $ab$ such that the minimum-area
triangles charged to $ab$ are fat. Note that the length of each edge
in $E$ is in the interval $[2/3^{1/4}, 2]$.

We next claim that, for any pair of points $p,q\in S$ with
$|pq| < 1$, neither $p$ nor $q$ can be an endpoint of an edge in $E$.
Indeed, suppose to the contrary that $p,q$ is such a pair and that $pa$
is an edge of $E$, for some $a\in S$; by construction, $a\ne q$.
Let $\Delta pab$ be a fat minimum-area triangle charged to $pa$.
If $q$ is collinear with $pa$, then $\Delta pqb$ is a nondegenerate
triangle of area strictly smaller than that of $\Delta pab$ (recall that
$|pq|< 1 < |pa|$), a contradiction.
If $q$ is not collinear with $pa$, $\Delta paq$ is a nondegenerate
triangle of area $ \leq \frac{|pa| \cdot |pq|}{2} < \frac{2 \cdot 1}{2} =1$,
again a contradiction.

Let $S' \subseteq S$ be the set obtained by repeatedly removing
the points of $S$ whose nearest neighbor in $S$ is at distance
smaller than $1$.  Clearly, the minimum inter-point distance in $S'$ is at
least $1$, and the endpoints of each edge in $E$ lie in $S'$. This implies,
via an easy packing argument, that the number of edges of $E$ incident
to any fixed point in $S'$ (all of length at most $2$) is only $O(1)$.
Hence $|E|=O(n)$. Since each edge in $E$ determines at most $10$
minimum-area triangles, as shown in the first part of our proof,
we conclude that $N_2=O(n)$, as claimed.

Hence there are at most $2{n\choose 2} + O(n) = n^2+O(n)$
minimum-area triangles in total.
\end{proof}

\section{Maximum-area triangles in 3-space}
\label{sec:maximum3}

\'Abrego and Fern\'andez-Merchant~\cite{af-02} showed that one can
place $n$ points on the unit sphere in $\RR^3$ so that they determine
$\Omega(n^{4/3})$ pairwise distances of $\sqrt{2}$
(see also~\cite[p.~191]{pa-95} and~\cite[p.~261]{bmp-05}).
This implies the following result:
\begin{theorem}\label{thm:maximum3}
For any integer $n$, there exists an $n$-element point set in $\RR^3$
that spans $\Omega(n^{4/3})$ triangles of maximum area, all incident to
a common point.
\end{theorem}
\begin{proof}
Denote the origin by $o$, and consider a unit sphere centered at
$o$. The construction in~\cite{af-02} consists of a set
$S=\{o\}\cup S_1\cup S_2$ of $n$
points, where $S_1\cup S_2$ lies on the unit sphere,
$|S_1|=\lfloor (n-1)/2\rfloor$, $|S_2|=\lceil (n-1)/2 \rceil$,
and there are $\Omega(n^{4/3})$ pairs of orthogonal segments of the form
$(os_i,os_j)$ with  $s_i \in S_1$ and $s_j \in S_2$.

Moreover, this construction can be realized in such a way that
$S_1$ lies in a small neighborhood of $(1,0,0)$, and $S_2$ lies in a
small neighborhood of $(0,1,0)$, say. The area of every right-angled
isosceles triangle $\Delta os_is_j$ with $s_i\in S_1$ and $s_j\in
S_2$ is $1/2$. All other triangles have
smaller area: this is clear if at least two vertices of a triangle
are from $S_1$ or from $S_2$; otherwise the area is given by
$\frac{1}{2}\sin \alpha$, where $\alpha$ is the angle of the two
sides incident to the origin, so the area is less than $1/2$ if
these sides are not orthogonal.
\end{proof}

\smallskip
We next show that the construction in Theorem~\ref{thm:maximum3} is
almost tight, in the sense that at most $O(n^{4/3+\eps})$ maximum-area
triangles can be incident to any point of an $n$-element point
set in $\RR^3$, for any $\eps>0$.

\begin{theorem}  \label{thm:max-3d}
  The number of triangles of maximum area spanned by a set $S$ of $n$
  points in $\RR^3$ and incident to a fixed point $a\in S$ is
  $O(n^{4/3+\eps})$, for any $\eps>0$.
\end{theorem}

Assume, without loss of generality, that the maximum area is $1$.
Similarly to the proof of Theorem~\ref{thm:unit3}, we map
maximum-area triangles to point-cylinder incidences. Specifically,
if $\Delta{abc}$ is a maximum-area triangle spanned by a point set $S$,
then every point of $S$ lies on, or in the interior of, the cylinder
with axis $ab$ and radius $2/|ab|$ ($c$ itself lies on the cylinder).
The following two lemmas give upper bounds on the number of
point-cylinder incidences in this setting. First we prove a
weaker bound (Lemma~\ref{lem:kst-cylinder}) which, combined with the
partition technique, gives an almost tight bound
(Lemma~\ref{lem:st-cylinder}).
%
Our proof is somewhat reminiscent of an argument of Edelsbrunner and
Sharir~\cite{es-91}, where it is shown that the number of point-sphere
incidences between $n$ points and $m$ spheres in $\RR^3$ is
$O(n^{2/3}m^{2/3}+n+m)$, provided that no point lies in the exterior
of any sphere.

\begin{lemma}\label{lem:kst-cylinder}
  Let $S$ be a set of $n$ points, and $\C$ a set of $m$ cylinders in
  $\RR^3$, such that the axis of each cylinder passes through the
  origin, and no point lies in the exterior of any cylinder. Then
  the number of point-cylinder incidences is
  $O(nm^{\frac{1+\eps}{2}}+m)$, for any $\eps>0$.
\end{lemma}
\begin{proof}
  Assume, without loss of generality, that the horizontal plane $h$
  incident to the origin does not contain any point of $S$, and that
  the points above $h$ participate in at least half of the
  point-cylinder incidences. For simplicity, continue to denote by
  $S$ the subset of the at most $n$ points lying above $h$.
  Consider the 3-dimensional dual arrangement $(S^*,\C^*)$, where
  the dual of a point $p\in \RR^3\setminus \{o\}$ is the cylinder
  $p^*$ with axis $op$ and radius $2/|op|$; and the dual of a cylinder
  $\gamma$ whose axis passes through the origin is a point $\gamma^*$
  above $h$ that lies on the axis of $\gamma$ at distance
  $2/{\rm radius}(\gamma)$ from the origin. Note that incidences
  between points and cylinders are preserved, and that no point of
  $\C^*$ lies in the exterior of any cylinder of $S^*$. It therefore
  suffices to prove that the number of incidences between $S^*$ and 
  $\C^*$ is
  $I(\C^*,S^*)=O(nm^{\frac{1+\eps}{2}}+m)$.

  Consider the intersection $B$ of the interiors of all cylinders in
  $S^*$. Since the interior of each cylinder is convex, $B$ is a convex
  body homeomorphic to a ball, whose boundary is composed of patches of
  cylinders. Faces, edges, and vertices of $B$ can be defined as
  connected components of the intersections of one, two, and three
  cylinders, respectively.
Each of the points of $\C^*$ that lie on faces of $\bd B$ contributes
one incidence. Since all the cylinder axes pass through the origin, no
edge of $\bd B$ can be straight, so it cannot be contained in any
cylinder of $S^*$ other than the two defining it (recall
Lemma~\ref{3int}). Hence the points of $\C^*$ that lie on faces
or edges of $\bd B$ contribute at most $2m$ incidences.

We are left with the task of bounding the
number of point-cylinder incidences involving points at vertices 
of $B$. Note that there may exist cylinders incident to a vertex $p$ 
of $B$ and not containing any other points of $\bd B$ in the 
vicinity of $p$. To account for such cylinders too, perturb
the radii of each cylinder in $S^*$, so that each radius $r$ is 
decreased to the radius $(1-\delta r)r$,
for a sufficiently small $\delta>0$ (that is, the radii of larger 
cylinders decrease by a larger factor).
As a result, every cylinder incident to a vertex $p\in \bd B$ is
replaced by a cylinder that defines a face in a sufficiently
small neighborhood of $p$ (even though it is not incident to
$p$ after this perturbation). The number of point-cylinder incidences 
between $\C^*$ and the vertices of $\bd B$ is proportional to
the number of vertices of the resulting $\bd B'$ after the perturbation. 
  By a result of Halperin and Sharir~\cite{hs-95}, the
  complexity of a single cell in the arrangement of $n$ constant
  degree algebraic surfaces in $\RR^3$ is $O(n^{2+\eps})$,
  for any $\eps>0$.
  Hence, we obtain an upper bound of
  $I(S,\C)=O(m+n^{2+\eps})$, for any $\eps>0$.

  Partition $S$ into $\lceil n/\sqrt{m}\rceil$ subsets, each
  containing at most $\sqrt{m}$ points. The preceding argument
  implies that each subset $S'\subset S$ has at most
  $I(S',\C)=O(m+(\sqrt{m})^{2+\eps}) = O(m^{1+\eps/2})$ incidences with
  the cylinders of $\C$. Therefore, altogether there are at most
  $\lceil n/\sqrt{m}\rceil \cdot O(m^{1+\eps/2})
  =O(nm^{\frac{1+\eps}{2}}+m)$ incidences.
\end{proof}

\begin{lemma}\label{lem:st-cylinder}
  Let $S$ and $\C$ be as in the preceding lemma.
  Then the number of point-cylinder
  incidences is $O((n^{2/3}m^{2/3}+n+m)^{1+\eps})$, for any $\eps>0$.
\end{lemma}
\begin{proof}
If $m>n^2$, then Lemma~\ref{lem:kst-cylinder} gives an upper bound
of $O(nm^{\frac{1+\eps}{2}}+m) = O(m^{1+\eps})$.
We may therefore assume henceforth that $m\leq n^2$.

  For an integer $r\in \NN$, to be specified later, choose a random
  sample $R\subset \C$ of $r$ cylinders, and let $B$ denote the
  intersection of the interiors of the cylinders in $R$.
  By~\cite{hs-95}, the combinatorial complexity of $B$ is
  $O(r^{2+\eps})$, for any $\eps>0$. Hence, the convex body $B$ can
  be partitioned into $O(r^{2+\eps})$ cells, each bounded by a
  constant number of constant-degree algebraic surfaces.
  (This can be done, e.g., by first partitioning $\bd B$ into
  pseudo-trapezoidal cells, and then by
  taking the convex hull of each cell on $\bd B$ with the origin.)
  By the $\eps$-net theory (see, e.g.,~\cite[Chap.~10.3]{m-02}), with
  constant probability, the interior of each cell intersects at most
  $O(m \frac{\log r}{r}) = O(m/r^{1-\eps})$ cylinders of $\C$.
  We may assume then that our sample $R$ has this property.
  Similarly to the proof of
  Lemma~\ref{lem:cut}, assign each point to a unique cell. Assign
  every point in the interior of a cell $\sigma_i$ to $\sigma_i$;
  assign every point on the boundary of several cells to the cell with
  minimum index. Let $n_i$ denote the number of points assigned to
  cell $\sigma_i$.

  Applying Lemma~\ref{lem:kst-cylinder} in each cell $\sigma_i$, we
  get the upper bound
  ${\displaystyle O\left(n_i
    \left(\frac{m}{r^{1-\eps}}\right)^{\frac{1+\eps}{2}} +
    \left(\frac{m}{r^{1-\eps}}\right)\right) }$
on the number of incidences between points assigned
  to $\sigma_i$ and cylinders intersecting the interior of $\sigma_i$.
  Summing over all $O(r^{2+\eps})$ cells, we have
  $$
  \sum_{i} O\left(n_i
    \left(\frac{m}{r^{1-\eps}}\right)^{\frac{1+\eps}{2}} +
    \left(\frac{m}{r^{1-\eps}}\right)\right) = O\left(n\left(
      \frac{m}{r^{1-\eps}}\right)^{\frac{1+\eps}{2}} +
    mr^{1+2\eps}\right) =
  O\left(\frac{nm^{\frac{1+\eps}{2}}}{r^{\frac{1-\eps}{2}}} +
    mr^{1+2\eps}\right)
  $$
  incidences of this kind. By choosing
  $r=\min\left\{\lfloor n^{2/3}/m^{1/3}\rfloor,\,m\right\}$, 
  this is at most
  $O(n^{2/3+\eps'}m^{2/3+\eps'}+n^{1+\eps'})$, for another, still
arbitrarily small, $\eps'>0$.  Finally, the number of incidences
  between points assigned to one cell and cylinders that do not
  intersect the interior of that cell can be bounded similarly to the
  proof of Lemma~\ref{lem:cut}: This number is proportional to the
  number of cells plus the number of points, which is
  $O(n+r^{2+\eps}) = O(n^{1+\eps})$, as is easily checked.
(In this final argument, we use the fact all axes pass through the
origin, so no 1-dimensional edge of $\bd B$ can be contained in more
than two cylinders; see also the proof of Lemma~\ref{lem:kst-cylinder}.)
\end{proof}

The upper bound of Lemma~\ref{lem:st-cylinder} is almost tight: For any 
$n$ and $m$, there are $n$ points and $m$ cylinders with axes through 
the origin and containing no points in their exterior, which determine 
$\Omega(n^{2/3}m^{2/3}+n+m)$ point-cylinder incidences. To construct
such a configuration, take $n$ points and $m$ lines on the plane 
$\pi:z=1$ in $\RR^3$ with $\Omega(n^{2/3}m^{2/3}+n+m)$ 
point-line incidences~\cite{st-83}. Project these points and lines centrally
from the origin onto the unit sphere, to obtain a system of $n$ points 
and $m$ great circles with the same number of incidences.
Each great circle of the unit sphere lies in a unique cylinder of unit 
radius whose axis passes through the origin, and every such cylinder 
contains all the other points of the unit sphere in its interior. 
This gives $n$ points on the unit sphere and $m$ cylinders of 
unit radius whose axes pass through the origin
(so that no point lies in the exterior of any cylinder),
with $\Omega(n^{2/3}m^{2/3}+n+m)$ point-cylinder incidences.

\begin{proofof}{Theorem~\ref{thm:max-3d}}
Let $A$ denote the maximum triangle area determined by a set $S$ of
$n$ points in $\RR^3$. For every point $a\in S$, consider the system
of $n-1$ points in $S\setminus \{a\}$ and $n-1$ cylinders, each defined
by a point $b\in S\setminus \{a\}$, and has axis $ab$ and radius
$2A/|ab|$.
Every point-cylinder incidence corresponds to a triangle of area $A$
spanned by $S$ and incident to $a$. Since $A$ is the maximum area, no
point of $S$ may lie in the exterior of any cylinder.  By
Lemma~\ref{lem:st-cylinder}, the number of such triangles is
$O(n^{4/3+\eps})$, for any $\eps>0$.
\end{proofof}

Theorems \ref{thm:maximum3} and \ref{thm:max-3d} imply the
following bounds on the number of maximum-area triangles in $\RR^3$:

\begin{theorem}\label{thm:allmax-3d}
The number of triangles of maximum area spanned by $n$ points in
$\RR^3$ is $O(n^{7/3+\eps})$, for any $\eps>0$.
For all $n \geq 3$, there exist $n$-element point sets in $\RR^3$
that span $\Omega(n^{4/3})$ triangles of maximum area.
\end{theorem}

\section{Distinct triangle areas in 3-space}\label{sec:dist3d}

Following earlier work by Erd\H{o}s and Purdy~\cite{ep-76}, Burton and
Purdy~\cite{bp-79}, and Dumitrescu and T\'oth~\cite{dt-07a},
Pinchasi~\cite{p-07} has recently proved that $n$ noncollinear points in
the plane always determine at least $\left\lfloor
  \frac{n-1}{2}\right\rfloor$ distinct triangle areas, which is
attained by $n$ equally spaced points distributed evenly on two
parallel lines. No linear lower bound is known in 3-space, and
the best we can show is the following:

\begin{theorem}\label{thm:dist3d}
  Any set $S$ of $n$ points in $\RR^3$, not all on a line, determines at
  least $\Omega(n^{2/3}/\beta(n))$ triangles of distinct areas, for
  some extremely slowly growing function $\beta(n)$.
  Moreover, all these triangles share a common side.
\end{theorem}

For the proof, we first derive a new upper bound
(Lemma~\ref{lem:clarkson}) on the number of point-cylinder
incidences in $\RR^3$, for the special case where the axes of the
cylinders pass through the origin (but without the additional requirement
that no point lies outside any cylinder).
Consider a set $\C$ of $m$ such cylinders.
These cylinders have only three degrees of
freedom, and we can dualize them to points in 3-space. We use a
duality similar to that used in the proof of Lemma~\ref{lem:kst-cylinder}.
Specifically, we fix some generic halfspace $H$ whose bounding plane
passes through the origin, say, the halfspace $z>0$. We then map
each cylinder with axis $\ell$ and radius $\varrho$ to the point on
$\ell\cap H$ at distance $1/\varrho$ from the origin;
and we map each point $p\in H$ to the cylinder whose axis is the line
spanned by $op$ and whose radius is $1/|op|$.
As argued above, this duality preserves point-cylinder incidences.

By (a dual version of) Lemma~\ref{3int}, any three points
can be mutually incident to at most eight cylinders whose axes
pass through the origin. That is, the
bipartite incidence graph (whose two classes of vertices correspond to
the points of $S$ and the cylinders of $\C$, respectively, and an edge
represents a point-cylinder incidence) is $K_{3,9}$-free. It follows
from the theorem of K\H{o}v\'ari, S\'os and Tur\'an~\cite{kst-54} (see
also~\cite[p.~121]{pa-95})
that the number of point-cylinder incidences is $O(nm^{2/3}+m)$.
We then combine this bound with the partition technique of
Clarkson~\etal~\cite{ceg-90}, to prove a sharper upper bound on
the number of point-cylinder incidences of this kind. Specifically, we
have:

\begin{lemma}\label{lem:clarkson}
  Given $n$ points and $m$ cylinders, whose axes pass through the
  origin, in 3-space, the number of point-cylinder incidences is
  $O(n^{3/4}m^{3/4}\beta(n)+n+m).$
\end{lemma}
\begin{proof}
Let $\C$ be the set of the $m$ given cylinders, and $S$ be the set of
the $n$ given points.
Let $h$ be a plane containing the origin, but no point of $S$, and
assume, without loss of generality, that the subset $S'$ of points
lying in the positive hafspace $h^+$ contributes at least half of
the incidences with $\C$.  If $m>n^3$, then the
K\H{o}v\'ari-S\'os-Tur\'an Theorem yields an upper bound of
$I(S',\C)=O(nm^{2/3}+m)=O(m)$. Similarly, if $m<n^{1/3}$,
the duality mentioned above leads to the bound
$I(S',\C)=O(mn^{2/3}+n)=O(n)$. For these two cases we have
then $I(S,\C) \leq 2I(S',\C)=O(m+n)$.
Assume henceforth that $n^{1/3}\leq m\leq n^3$.

We apply Lemma~\ref{lem:cut} with parameter 
$r= \lfloor n^{3/8}/m^{1/8} \rfloor$, and
use the K\H{o}v\'ari-S\'os-Tur\'an Theorem to bound the number of
incidences between the at most $n/r^3$ points and $m/r$ cylinders
in each subproblem.
Note that $1 \leq r \leq m$ in the above range of $m$.
The total number of incidences is thus
\begin{eqnarray*}
I(S,\C) &=& O(n+mr^2\beta(r)) + O(r^3\beta(r))\cdot
O \left( \frac{n}{r^3} \cdot \left(\frac{m}{r}\right)^{2/3}
+\frac{m}{r} \right)\\
&=&
O \left( n+ \frac{m^{2/3}n}{r^{2/3}}\beta(n) + m r^2\beta(r)\right) =
O\left(n+ n^{3/4}m^{3/4}\beta(n)\right).
\end{eqnarray*}
Putting all three cases together gives the bound in the theorem.
\end{proof}

\begin{proofof}{Theorem~\ref{thm:dist3d}}
  If there are $n/100$ points in a plane but not all on a line, then
  the points in this plane already determine $\Omega(n)$ triangles of
  distinct areas~\cite{bp-79}. We thus assume, in the remainder of the
  proof, that there are at most $n/100$ points on any plane.

  According to a result of Beck~\cite{b-83}, there is an absolute
  constant $k\in \NN$ such that if no line is incident to $n/100$ points
  of $S$, then $S$ spans $\Theta(n^2)$ distinct lines, each of which
  is incident to at most $k$ points of $S$. Since each point of $S$ is
  incident to at most $n-1$ of these lines, there is a point $a\in S$
  incident to $\Theta(n)$ such lines. Select a point of $S\setminus \{a\}$
  on each of these lines, to obtain a set $P$ of $\Theta(n)$ points.

  Let $t$ denote the number of distinct triangle areas determined by
  $S$, and let $\alpha_1,\alpha_2,\ldots ,\alpha_t$ denote these
  areas.  For each point $b\in P$ and $i=1,2,\ldots ,t$, we define a
  cylinder $C(ab,\alpha_i)$ with axis (the line spanned by) $ab$ and radius
  $2\alpha_i/|ab|$.
  Every point $c\in S$ for which the area of the
  triangle $\Delta{abc}$ is $\alpha_i$ must lie on the cylinder
  $C(ab,\alpha_i)$. Let $\C$ denote the set of the $O(nt)$ cylinders
  $C(ab,\alpha_i)$, for $b\in P$ and $i=1,2,\ldots,t$.
  For each point $b\in P$, there are
  $n-k=\Theta(n)$ points off the line through $ab$, each of which must
  lie on a cylinder $C(ab,\alpha_i)$ for some $i=1,2,\ldots , t$.
  Therefore, the number $I(S,\C)$ of point-cylinder incidences between
  $S$ and $\C$ is $\Omega(n^2)$. On the other hand, by
  Lemma~\ref{lem:clarkson}, we have
$$\Omega(n^2)\leq I(S,\C)\leq O(n^{3/4}(nt)^{3/4}\beta(n)+n+nt)
= O(n^{3/2}t^{3/4}\beta(n)),$$
which gives $t=\Omega (n^{2/3}/\beta^{4/3}(n))=\Omega
(n^{2/3}/\beta'(n))$, for another function $\beta'(n)$ of the same
slowly growing type, as required.
\end{proofof}

\section{Conclusion} \label{sec:conclusion}

We have presented many results on the number of triangles of specific areas
determined by $n$ points in the plane or in three dimensions. Our results improve
upon the previous bounds, but, most likely, many of them are not
asymptotically tight. This leaves many open problems of closing the
respective gaps. 
Even in cases where the bounds are asymptotically tight, such as
those involving minimum-area triangles in two and three dimensions,
determining the correct constants of proportionality still offers
challenges.

Here is yet another problem on triangle areas, of a slightly
different kind, with triangles determined by lines, not points
(motivated in fact by the question of bounding $|U_2|$ in the
proof of Theorem \ref{thm:unit2}). Any three nonconcurrent, and
pairwise non-parallel lines in the plane determine a triangle of
positive area. What is the maximum number of unit area
triangles determined by $n$ lines in the plane?

\begin{theorem}\label{thm:lines}
The maximum number of unit-area triangles determined by $n$ lines in 
the plane is $O(n^{7/3})$, and for any $n \geq 3$, there are $n$ lines that
determine $\Omega(n^2)$ unit-area triangles. 
\end{theorem}
\begin{proof}
{\em Lower bound}: Place $n/3$ equidistant parallel lines at angles 
$0$, $\pi/3$, and $2\pi/3$, through the points of an appropriate section of 
the triangular lattice, and observe that
there are $\Omega(n^2)$ equilateral triangles of unit side (i.e., of
the same area) in this construction.

{\em Upper bound}: Let $L$ be a set of $n$ lines in the plane. 
We define a variant
of the hyperbolas used in the proof of Theorem \ref{thm:unit2}:
For any pair of non-parallel lines $\ell_1,\ell_2\in L$, let
$\gamma(\ell_1,\ell_2)$ denote the locus of points $p\in \RR^2$,
$p\not\in\ell_1\cup \ell_2$, such that the parallelogram that has a
vertex at $p$ and two sides along $\ell_1$ and $\ell_2$,
respectively, has area $1/2$. The set $\gamma(\ell_1,\ell_2)$ is the
union of two hyperbolas with $\ell_1$ and $\ell_2$ as asymptotes 
(four connected branches in total).
Any two non-parallel lines uniquely determine two such hyperbolas.
Let $\Gamma$ denote the set of the branches of these hyperbolas, and
note that $|\Gamma|=O(n^2)$.
Observe now that, if $\ell_1$, $\ell_2$, and $\ell_3$ determine a
unit area triangle, then $\ell_3$ is tangent to one of the two
hyperbolas in $\gamma(\ell_1,\ell_2)$.

We first derive a weaker bound. Construct two bipartite graphs 
$G_1,G_2\subseteq L\times\Gamma$. We put an edge $(\ell,\gamma)$ 
in $G_1$ (resp., $G_2$) if $\ell$ is tangent to $\gamma$ and 
$\ell$ lies below (resp., above) $\gamma$. The edges of $G_1$ 
and $G_2$ account for all line-curve tangencies. 
Observe that neither graph contains a $K_{5,2}$, 
that is, there cannot be five distinct lines in $L$ tangent 
to two branches of hyperbolas from above (or from below). 
Indeed, this would force the two branches to intersect at five
points, which is impossible for a pair of distinct quadrics.
It thus follows from the theorem of 
K\H{o}v\'ari, S\'os and Tur\'an~\cite{kst-54} (see 
also~\cite[p.~121]{pa-95}) that the number of line-hyperbola tangencies 
between any $n_0$ lines in $L$ and any $m_0$ hyperbolas in $\Gamma$ 
is $O(n_0 m_0^{4/5}+m_0)$. With $n_0=n$ and $m_0=O(n^2)$, this already
gives a bound of $O(n \cdot n^{8/5}+n^2)=O(n^{13/5})$ on the number of
unit-area triangles determined by $n$ lines in the plane. We next
derive an improved bound.  

Let $L$ be the given set of $n$ lines, and let $\Gamma$ be the
corresponding set of $m=O(n^2)$ hyperbola branches. 
We can assume that no line in $L$ is vertical, and apply a standard
duality which maps each line $\ell\in L$ to a point $\ell^*$. A
hyperbolic branch $\gamma$ is then mapped to a curve $\gamma^*$,
which is the locus of all points dual to lines tangent to $\gamma$;
it is easily checked that each $\gamma^*$ is a quadric.
Let $L^*$ denote the set of the $n$ dual points, and let
$\Gamma^*$ denote the set of $m=O(n^2)$ dual curves. 
A line-hyperbola tangency in the primal plane is then
mapped to a point-curve incidence in the dual plane. 

We next construct a $(1/r)$-cutting for $\Gamma^*$, partitioning
the plane into $O(r^2)$ relatively open cells of bounded description
complexity, each of which contains at most $n/r^2$ points and is
crossed by at most $m/r$ curves. By using the previous bound for each
cell, the total number of incidences involving points in the interior
of these cells is
$$ O\left( 
r^2 \left(\frac{n}{r^2} \left(\frac{m}{r}\right)^{4/5} + \frac{m}{r}\right) 
\right)= O\left(n \left(\frac{m}{r}\right)^{4/5} + mr\right).
$$ 
We balance the two terms by setting $r=n^{5/9}/m^{1/9}$, and observe
that $1\le r\le m$ if $m\le n^5$ and $n\le m^2$; since
$m=\Theta(n^2)$, both inequalities do hold in our case. Hence, the
total number of incidences under consideration is
$O(m^{8/9}n^{5/9})=O(n^{7/3})$.

It remains to bound the overall number of incidences involving points 
lying on the boundaries of at least two cells. A standard
argument, which we omit, shows that the number of these incidences
is also $O(n^{7/3})$, and thereby completes the proof of the theorem.
\end{proof}

Some remarks are in order: The line variant of unit-area triangle 
problems is {\em not} equivalent to the point variant, under the 
standard point-line duality.  Specifically: Let $S$ be a set of $n$ 
points in the plane having distinct $x$-coordinates. Consider the 
duality transform that maps a point
$p=(a,b)$ to the line $p^*:\, y=ax-b$, and vice versa.
It is easy to see that
there is no absolute constant $A>0$ such that, for $p,q,r \in S$,
triangle $\Delta{pqr}$ has unit area if and only if the triangle
$\Delta {p^*q^*r^*}$ formed by the three dual lines has area $A$.

Yet, there is a connection between the point- and the line-variants of
the unit-area problem in the plane. 
Go back to the notation in the proof of Theorem \ref{thm:unit2},
where, for a parameter $k \leq n^{1/3}$, we had $|U_1|= O(n^2 k)$.
Recall that $U_2$ denotes the set of unit-area triangles where
all three top lines are {\em $k$-rich}, and that there are 
$|L_k|=O(n^2/k^3)$ such lines. Observe that the three top lines 
of each triangle in $U_2$ determine a triangle of area $4$. We thus 
face the question of bounding the number of triangles of area 4 
determined by the $k$-rich lines in $L_k$. 
By Theorem~\ref{thm:lines}, there are most $O((n^{2}/k^{3})^{7/3})$
such triangles. Balancing $|U_1|$ with $|U_2|$ yields
$k=n^{1/3}$, thereby implying that $|U_1|+|U_2| = O(n^{7/3})$.

We note that the bound $O(n^{44/19})$ of Theorem \ref{thm:unit2}
could be re-derived with this new approach, if the bound of
Theorem~\ref{thm:lines} could be improved to $O(n^{11/5})$.
Moreover, an $o(n^{11/5})$ bound for the line-variant would in turn
lead to an improvement in our current bound for the classical
point-variant of the unit area problem in the plane.

\newpage
\section*{Appendix}

\begin{proofof}{Lemma~\ref{3int}}
Let us recall from~\cite{EPS:cyl} the structure of the
intersection curve between two cylinders.  Let $C$ and $C'$ be two
cylinders with nonparallel axes, so each pair of axes are either skew
to each other or concurrent.  Let $\gamma$ denote the curve of their
intersection.

To simplify the analysis, we assume, without loss of generality, that
the axis $\alpha$ of $C$ is the $z$-axis and that its radius is $1$.
Let $\alpha'$ and $\rho'$ denote respectively the axis and radius of
$C'$. Let $\pi$ be the plane passing through $\alpha'$ and through the
shortest segment $e$ connecting the axes $\alpha,\alpha'$. If
$\alpha,\alpha'$ are skew lines, $e$ and $\pi$ are well defined.  If
$\alpha$ and $\alpha'$ are concurrent, we take $\pi$ to be the plane
passing through $\alpha'$ and orthogonal to the plane spanned by
$\alpha$ and $\alpha'$.

Let $\sigma$ denote the ellipse $C\cap\pi$. We use a cylindrical
coordinate system $\theta,z$ on $C$, and write the equation of
$\sigma$ as $z=a\cos\theta+b\sin\theta+c$, where $z=ax+by+c$ is
the quation of $\pi$.

As shown in \cite{EPS:cyl}, the equation of $\gamma$ is
$$
z=\sigma(\theta)\pm
  \frac{1}{\sin\beta}\sqrt{(\rho')^2-d^2(\sigma(\theta),\alpha')} ,
$$
where $\beta$ is the angle between the axes. Moreover,
$d(\sigma(\theta),\alpha')$, being the distance, within $\pi$,
of a point on the ellipse $\sigma$ from the line $\alpha'$,
can also be expressed as $|p\cos\theta+q\sin\theta+r|$, for
appropriate parameters $p,q,r$.

Let now $C,C_1,C_2$ be three cylinders with no pair of parallel axes.
Suppose to the contrary that $|C\cap C_1\cap C_2|\ge 9$.
Let $\gamma_i$ denote the intersection curve $C\cap C_i$, for
$i=1,2$. Write the equations of $\gamma_1,\gamma_2$ as
$$
z=a_i\cos\theta+b_i\sin\theta+c_i \pm
  \frac{1}{\sin\beta_i}\sqrt{(\rho_i)^2-(p_i\cos\theta+q_i\sin\theta+r_i)^2},
$$
for $i=1,2$, with the appropriate parameters as above.
We can re-parameterize these curves by putting $t=\tan(\theta/2)$
and $w=z(1+t^2)$, to obtain two equations of the form
\begin{eqnarray*}
w & = & Q_1(t)\pm\sqrt{K_1(t)} \\
w & = & Q_2(t)\pm\sqrt{K_2(t)} ,
\end{eqnarray*}
where $Q_1,Q_2$ are quadratic polynomials and $K_1,K_2$ are
quartic polynomials. We are given that these two equations have
at least $9$ common roots (it is easy to check that distinct roots of
the original system are mapped to distinct roots of the new system).

If $Q_1(t)\equiv Q_2(t)$ then the common roots must satisfy
$K_1(t)=K_2(t)$. Since there are at least $9$ such roots and this
is a quartic equation, we must also have $K_1(t)\equiv K_2(t)$.

We will get to this case soon, but let us first consider the case
$Q_1(t)\not\equiv Q_2(t)$. After squaring, the
equations become
\begin{eqnarray*}
(w-Q_1(t))^2 & = & K_1(t) \\
(w-Q_2(t))^2 & = & K_2(t) .
\end{eqnarray*}
Hence
$$
w = -\frac{K_2(t)-K_1(t)}{2(Q_2(t)-Q_1(t))} + \frac{Q_1(t)+Q_2(t)}{2} ,
$$
so $t$ must satisfy the equation
\begin{equation} \label{deg8eq}
\left(
 -\frac{K_2(t)-K_1(t)}{2(Q_2(t)-Q_1(t))} + \frac{Q_2(t)-Q_1(t)}{2}
\right)^2 = K_1(t) ,
\end{equation}
which is a polynomial equation of degree at most $8$. Since it has $9$
roots, it must vanish identically.

Since the left-hand side of (\ref{deg8eq}) is a square, $K_1$ must
also be a square.  However, $K_1(t)$ is proportional to
$$
\biggl(\rho_1(1+t^2)\biggr)^2 -
\biggl(p_1(1-t^2)+2q_1t+r_1(1+t^2)\biggr)^2 =
$$
$$
\biggl(\rho_1(1+t^2)-(p_1(1-t^2)+2q_1t+r_1(1+t^2))\biggr) \cdot
\biggl(\rho_1(1+t^2)+(p_1(1-t^2)+2q_1t+r_1(1+t^2))\biggr) .
$$
It follows that either each of these factors is a square, or they
are multiples of each other. In the former case, we must have
\begin{eqnarray*}
q_1^2 & = & (\rho_1+p_1-r_1)(\rho_1-p_1-r_1) = (\rho_1-r_1)^2-p_1^2 \\
q_1^2 & = & (\rho_1-p_1+r_1)(\rho_1+p_1+r_1) = (\rho_1+r_1)^2-p_1^2 ,
\end{eqnarray*}
implying that $\rho_1-r_1 = \pm (\rho_1+r_1)$, so either $\rho_1=0$
or $r_1=0$. The first equality is impossible---our cylinders have
positive radii. The second equality implies that
$\rho_1^2=p_1^2+q_1^2$. However, as argued in \cite{EPS:cyl}, by
shifting $\theta$, we may assume that $q_1=0$ and $p_1$ is half the
major axis of $\sigma_1$. This implies that $\sigma_1$ is a circle
(since its minor axis is always equal to $2\rho_1$),
which can happen only when $\alpha_1$ is orthogonal to $\alpha$.
Moreover, $r_1=0$ implies that $\alpha$ and $\alpha'$ are concurrent.

In the latter case, since $\rho_1\ne 0$, the two factors are
proportional to each other only when $p_1(1-t^2)+2q_1t$ is a
multiple of $1+t^2$, which can only happen when $p_1=q_1=0$,
which again is impossible.

Since the only remaining case is that of orthogonal concurrent axes,
it follows, using a
symmetric argument, that in the only remaining case, the three axes
$\alpha,\alpha_1,\alpha_2$ are concurrent, at a common point, and
mutually orthogonal. It is easily checked that in this case the cylinders
can intersect in at most $8$ points, contrary to assumption.
(This special case of three intersecting cylinders has been studied
a lot; see, e.g.,~\cite{Baumann}.)

Hence, $Q_1(t)\equiv Q_2(t)$ and $K_1(t)\equiv K_2(t)$. However, the
first identity implies that $\sigma_1=\sigma_2$,
so the plane containing the axis of $C_1$
also contains the axis of $C_2$. Since these axes are nonparallel,
they must be concurrent. Since the analysis is fully symmetric with
respect to the three cylinders, it follows that all three axes are
either coplanar or concurrent. If they are coplanar but not concurrent,
then it is easy to check that the planes $\pi_1$ and $\pi_2$ (with
respect to $C$ as the ``base'' cylinder) cannot be equal. If the
three axes are concurrent then again the identity of the planes
$\pi_1,\pi_2$ implies that both $\alpha_1$ and $\alpha_2$ must be
orthogonal to $\alpha$, and the fact that the argument is fully symmetric
implies that all three axes must be concurrent and mutually orthogonal,
a case that we have already ruled out.
This completes the proof of Lemma~\ref{3int}.
\end{proofof}

\end{document}